\documentclass[a4paper,11pt,twoside]{article}
\usepackage[left=2cm,right=2cm,top=2cm,bottom=3cm]{geometry}

\usepackage[T1]{fontenc}
\usepackage[utf8]{inputenc} 
\usepackage[english]{babel} 

\usepackage[super]{nth}  
\usepackage{units}
\usepackage{amsmath}

\usepackage{comment}
\usepackage{xcolor}

\definecolor{forestgreen}{rgb}{0.0, 0.5, 0.0}

\usepackage[usestackEOL]{stackengine}
\usepackage{amsmath} 
\usepackage{amssymb} 
\usepackage{amsthm}
\usepackage{amsfonts}
\usepackage{mathtools}
\usepackage{mathrsfs}
\usepackage{braket} 
\usepackage{graphicx}
\usepackage{graphics}
\usepackage{stmaryrd}
\usepackage{textcomp}
\usepackage[export]{adjustbox}
\usepackage{cancel}
\usepackage{mathabx,epsfig}
\usepackage{cases}
\usepackage[overload]{empheq}
\mathtoolsset{showonlyrefs=true,showmanualtags=true}

\usepackage[all]{xy}

\newcommand*{\Scale}[2][4]{\scalebox{#1}{$#2$}}

\newcommand{\numberset}{\mathbb}
\newcommand{\N}{\numberset{N}}

\newcommand{\R}{\numberset{R}}

\newcommand{\E}{\numberset{E}}
\DeclareMathOperator{\law}{\text{Law}}  
\DeclareMathOperator{\w}{\mathcal{W}}  
\DeclareMathOperator{\pr}{\mathscr{P}} 
\DeclareMathOperator{\p}{\numberset{P}}  
\DeclareMathOperator{\meas}{meas}   
\newcommand{\ndelta}{M}    
\newcommand{\ngrid}{N}  
\newcommand{\doublehookrightarrow}{\mathrel{\mathrlap{{\mspace{4mu}\lhook}}{\hookrightarrow}}}

\newcommand{\dirnum}{d_v}
\newcommand{\dirn}{d_v}
\newcommand{\lin}{\mathscr{L}}
\DeclareMathOperator{\tone}{\theta_1}
\DeclareMathOperator{\ttwo}{\theta_2}
\DeclareMathOperator{\kone}{\kappa_1}
\DeclareMathOperator{\ktwo}{\kappa_2}
\DeclareMathOperator{\e}{e_{\alpha}}
\DeclareMathOperator{\eplus}{e_{\alpha}^+}
\DeclareMathOperator{\h}{\mathscr{H}}
\DeclareMathOperator{\cmn}{c_{MN}}
\DeclareMathOperator{\esssup}{\text{ess-sup}}

\DeclareMathOperator{\cm}{c_{MN}}

\theoremstyle{definition}
\newtheorem{definition}{Definition}[section]
\newtheorem{theorem}[definition]{Theorem}
\newtheorem{lemma}[definition]{Lemma}

\newtheorem{proposition}[definition]{Proposition}

\newtheorem{remark}[definition]{Remark}

\newtheorem{setting}[definition]{Setting}
\newtheorem*{notation*}{Notation}
\newtheorem*{definition*}{Definition}
\newtheorem*{theorem*}{Theorem}
\newtheorem*{lemma*}{Lemma}
\newtheorem*{corollary*}{Corollary}
\newtheorem*{proposition*}{Proposition}
\newtheorem*{fact*}{Fact}
\newtheorem*{example*}{Example}
\newtheorem*{claim*}{Claim}
\newtheorem*{remark*}{Remark}
\newtheorem*{conjecture*}{Conjecture}

\numberwithin{equation}{section}

\allowdisplaybreaks[4]

\title{On the fluctuations of an SDE system modelling grid cells}
\author{Andrea Clini\thanks{Mathematical Institute, University of Oxford, Oxford OX2 6GG, UK (andrea.clini@maths.ox.ac.uk)}
}
\date{\today}

\allowdisplaybreaks[4]

\begin{document}

\maketitle

\begin{abstract}
Several differential equation models have been proposed to explain the formation of patterns characteristic of the grid cell network. Understanding the effect of noise on these models is one of the key open questions in computational neuroscience. 
In the present work, we continue the analysis of the SDE system \eqref{abstract model particle system} initiated in \cite{carrillo_clini_solem_the_mean_field_limit}.
We show that the fluctuations of the empirical measure associated to \eqref{abstract model particle system} around its mean field limit converge to the solution of a Langevin SPDE.
The interaction between different columns of neurons along the cortex prescribes a peculiar scaling regime.
\end{abstract}





\section{Introduction}
\label{introduction}

Grid cells are a particular type of neuron in the brain of mammals discovered in 2005 \cite{gridcells}, see also the review \cite{McNaughtonMoser}.
These neurons fire at regular intervals as an animal moves across an area, storing information such as position, direction and velocity, and thus enabling it to understand its movement in space.

Since their discovery, there has been extensive research to understand the precise behavior of grid cells, see \cite{tenyears,McNaughtonMoser,bressloff_spatiotemporal_dynamics_of_continuum_neural_fields} and the references therein.
Mathematically, their network is commonly described by deterministic continuous attractor network dynamics through a system of neural field models \cite{Ermentrout2010, McNaughtonetal, burakfiete, coueyetal}, based on the classical papers \cite{WC1, WC2, Amari1977}. 
In particular, as the brain is inherently noisy, understanding how the grid cell network is affected by noise is one of the currently open challenges in the field.
This question has recently been addressed from several directions \cite{bresslof_stochastic_neural_field_model_of_stimulus,BurakFieteNoise,BAC19,KE13,MB,TouboulPhysD,CHS2020,carrillo_roux_solem_noise_driven_bifurcation}.

The present work is the continuation of \cite{carrillo_clini_solem_the_mean_field_limit}, along the direction initiated in \cite{CHS2020}.
In \cite{CHS2020} the authors studied a system of Fokker--Planck PDEs derived by adding noise to the attractor network models in \cite{burakfiete, coueyetal} and formally taking the mean field limit.
The PDE has been further analyzed in \cite{carrillo_roux_solem_noise_driven_bifurcation}.
The passage to the mean field limit has been rigorously proved in \cite{carrillo_clini_solem_the_mean_field_limit}, where the authors derive the limit for a generalized model, covering diverse concrete models commonly proposed
(cf. \cite{AB20,burakfiete,BurakFieteNoise,coueyetal} and the review \cite{bressloff_spatiotemporal_dynamics_of_continuum_neural_fields}).
The aim of the present work is to analyze the fluctuations of the empirical measure of the generalized noisy network model around its deterministic mean field limit and prove a central limit theorem (cf. Theorem \ref{theorem/ clt/ clt section} below).

Precisely, we consider the following model for the interaction among $NM$ neurons stacked along $N$ columns, with $M$ neurons each, at the locations $x_1,\dots,x_N\in Q$ in a region $Q$ of the neural cortex.  
The $k^{th}$ neuron at location $x_i$ has an activity level $u_{ik}^\beta$ for each possible orientation $\beta=1,\dots,\dirnum$, corresponding to spatial directions (typically $\dirnum=4$).
The total activity level $u_{ik}=(u_{ik}^1,\dots,u_{ik}^{\dirnum})$ evolves according to the system
\begin{subequations}\label{abstract model particle system}
\begin{align}[left ={\empheqlbrace}]
    u_{ik}(t)=&\,u_{ik}(0)+\int_0^t\!b(x_i,r,u_{ik}(r),f_{N,\ndelta}(r))\,dr+\int_0^t\!\sigma(x_i,r,u_{ik}(r),f_{N,\ndelta}(r))\,dW_{ik}(r)-\ell_{ik}(t), \label{abstract model particle first line}
    \\[4mm]
    \ell_{ik}^{\beta}(t)=&\,-|\ell_{ik}^{\beta}|(t),\quad |\ell_{ik}^{\beta}|(t)=\int_0^t1_{\{u^\beta_{ik}(r)=0\}}d|\ell_{ik}^{\beta}|(r)\quad\text{for }\beta=1,\dots,\dirnum. \label{abstract model particle second line}
\end{align}
\end{subequations}
For simplicity, we take the cortex to be $Q=[0,1]^d$ and we extend everything periodically out of $Q$, since the network of neurons is commonly considered to have toroidal connectivity.
The results in this work are easily extended to any bounded open subset $Q\subseteq\R^d$, for any $d\geq1$.
Here, for integers $k=1,\dots,M$, we have i.i.d. families of random initial conditions $\{u_{ik}(0)\}_{i=1,\dots,N}$ for each space point $x_i$ in the cortex $Q$.
Moreover, for integers $i=1,\dots, N$ and $k=1,\dots, M$, we have $\dirnum$-dimensional Brownian motions $(W_{ik}^\beta)_{\beta=1,\dots,\dirnum}$,which can also be correlated.

The precise shape and properties of the coefficients $b$ and $\sigma$ are given in Section \ref{section/Hypotheses and notations} and are deduced from those of the concrete models considered (see the introduction of \cite{carrillo_clini_solem_the_mean_field_limit}).
We also denoted $f_{N,M}$ the empirical measure associated to these particles, that is
\begin{equation}\label{formula/ single empirical measure}
    f_{\ngrid,\ndelta}(r,dy,du)=\frac{1}{\ngrid\ndelta}\sum_{j=1}^{\ngrid}\sum_{m=1}^\ndelta\delta_{(x_j,u_{jm}(r))}\quad\text{regarded as a measure on $Q\times\R^{\dirnum}$.}
\end{equation}

Finally, for each $i$, $k$ and $\beta$, the term $\ell_{ik}^\beta$ is a finite variation process defined by \eqref{abstract model particle second line} which prevents the activity level $u_{ik}^\beta$ from taking negative values. 
Namely, as we can see in its definition, at each time $t$ this process equals the opposite of its total variation $\ell_{ik}^{\beta}(t)=-\big|\ell_{ik}^{\beta}\big|(t)$. In turn, the total variation stays constant when $u_{ik}^\beta>0$ and it increases in the form $\big|\ell_{ik}^{\beta}\big|(t)=\int_0^t1_{\{u_{ik}^\beta(r)=0\}}d\big|\ell_{ik}^{\beta}\big|(r)$ when $u_{ik}^\beta=0$, so as to push  $u_{ik}^\beta$ away from zero which is being dragged by the other terms on the right hand side of \eqref{abstract model particle first line}.
The introduction of such terms and constraints is therefore known as imposing \emph{reflecting boundary conditions} and $\ell_{ik}^\beta$ is called a \emph{reflection term}.
The existence and uniqueness of such a term need of course to be proved and this process is often referred to as the \emph{Skorokhod problem}.
Precise details concerning the well-posedness and the construction of the reflection term in our setting are all presented in the seminal papers \cite{Lions-Sznitman-1984-SDEreflectingBC, Sznitman1984} by Lions and Sznitman.

\begin{remark}
\label{remark/n mean field problem}
Rather than being a classical mean field system with $MN$ particles, it is more accurate to regard the model \eqref{abstract model particle system} as a collection of $N$ mean field systems with $M$ particles, one for each column of neurons, interacting among themselves.
Different neurons in the same column are subjected to independent sources of noise and are initiated with i.i.d. data.
However, different columns feature possibly correlated initial data and noise sources, so that standard independence arguments cannot be exploited to handle the column-column interaction and we have to rely on the continuity of the quantities involved.
This feature is well-analyzed in \cite{carrillo_clini_solem_the_mean_field_limit} and results in the convergence rate $\nicefrac{1}{\sqrt{M}}+\nicefrac{1}{N^{\nicefrac{\alpha}{d}}}$ for the mean field limit $M,N\to\infty$ (cf. Theorem \ref{theorem/Mean squared error estimates for actual particles vs McKean--Vlasov particles} below) in contrast to the usual decay rate $\nicefrac{1}{\sqrt{MN}}$ for a mean field problem with $MN$ particles.
In turn, this aspect will have consequences on the correct scaling regime needed to see a nontrivial behavior of the fluctuations (cf. Remark \ref{remark/ rescaling factor sqrt M} and Theorem \ref{theorem/ clt/ clt section} below).
\end{remark}

Under suitable assumptions, strong existence and uniqueness for \eqref{abstract model particle system} is proved in \cite[Theorem 2.3]{carrillo_clini_solem_the_mean_field_limit}.
In turn, the authors rigorously analyze the mean field limit of \eqref{abstract model particle system} as the number of columns $\ngrid$ and the number of neurons per column $\ndelta$ go to infinity.
It is shown (cf. Theorem \ref{theorem/Mean squared error estimates for actual particles vs McKean--Vlasov particles} below) that the particles converge to the solution of the family of McKean--Vlasov SDEs, for $x\in Q$,
\begin{subequations}\label{abstract 4 model McKean--Vlasov}
    \begin{align}[left ={\empheqlbrace}]
    \Bar{u}^{\epsilon}(x,t)=&u(x,0)+\int_0^tb(x,r,\Bar{u}^{\epsilon}(x,r),f(r))\,dr+\int_0^t\sigma(x,r,\Bar{u}^{\epsilon}(x,r),f(r))\,dW^{\epsilon}(x,r)-\bar{\ell}(x,t), \label{abstract 4 model McKean--Vlasov first line}
    \\
    \bar{\ell}^{\beta}(x,t)=&-|\bar{\ell}^{\beta}(x,\cdot)|(t),\quad |\bar{\ell}^{\beta}(x,\cdot)|(t)=\int_0^t1_{\{(\Bar{u}^{\epsilon})^\beta(x,r)=0\}}\,d|\bar{\ell}^{\beta}(x,\cdot)|(r)\quad\text{for }\beta=1,\dots,\dirnum. \label{abstract 4 model McKean--Vlasov second line}
    \end{align}
\end{subequations}
Here we have a family of random initial data $(u(x,0))_{x\in Q}$ for each point $x$ in the cortex.
The term $W^\epsilon(x,t)$ denotes a Gaussian noise term, white in time and $\epsilon$-correlated in space.
The space correlation is needed since equation \eqref{abstract 4 model McKean--Vlasov first line} would be ill-posed for an actual space-time white noise, and in fact it is also meaningful from the modelling point of view (cf. \cite[Remark 1.2]{carrillo_clini_solem_the_mean_field_limit} and Remark \ref{remark/ joint scaling regime in m n epsilon} below).
The precise definition and further properties of $W^{\epsilon}(x,t)$ are given in Section \ref{section/Hypotheses and notations}.

The coefficients $b$ and $\sigma$ are the same as in \eqref{abstract model particle system}, but here we have set $f(t,y,du)\coloneqq\law_{\R^{\dirnum}}(\Bar{u}(y,t))$, considered as a measure on $\R^{\dirnum}$ depending on $t\in[0,\infty)$ and $y\in Q$. 
This in turn induces a probability measure $f(t,dx,du)$ on $Q\times\R^{\dirnum}$ defined by integration as 
\begin{equation}\label{formula/ measure induced via integration}
    \int_{Q\times\R^{\dirnum}}\varphi(x,u) f(t,dx,du)\coloneqq\int_Q\int_{\R^{\dirnum}}\varphi(x,u)f(t,x,du)\,dx\qquad\text{for any $\varphi\in C_b(Q\times\R^{\dirnum})$}.
\end{equation}

Finally, for each fixed $x\in Q$ and $\beta=1,\dots,\dirnum$, the finite variation process $\bar{\ell}^\beta(x,t)$ defined in \eqref{abstract 4 model McKean--Vlasov second line} is again the reflection term coming from the Skorokhod problem (see the explanation after equation \eqref{abstract model particle system}) and it ensures that $(\bar{u}^{\epsilon})^\beta(x,t)\geq0$ for each $x$, $t$ and $\beta$. We refer the reader to \cite{Sznitman1984} for the details about such a process in the context of a classical McKean--Vlasov equation.

A formal application of It\^o formula shows that the law $f(t,x,du)\coloneqq\law_{\R^{\dirnum}}(\Bar{u}^\epsilon(x,t))$ is a weak solution of the Fokker--Planck equation with no-flux boundary conditions
\begin{subequations}
\label{abstract model nonlinear fokker--planck}
\begin{align}[left ={\empheqlbrace}]
&\partial_tf(t,x,u)\!+\!\!\nabla_u\!\cdot\!\Big(b(x,t,u,f(t))f(t,x,u)\Big)\!=\frac{1}{2}\sum_{\beta=1}^{\dirnum}\frac{\partial^2}{\partial u^\beta\partial u^\beta}\Big(\sigma_\beta(x,t,u,f(t))^2f(t,x,u)\Big), \label{abstract model nonlinear fokker--planck first line}
\\
&b_\beta(t,x,u,f(t))f(t,x,u)\!-\frac{1}{2}\frac{\partial}{\partial u^\beta}\Big(\sigma_{\beta}(x,t,u,f(t))^2f(t,x,u)\Big)~\Big|_{u^\beta=0}\!\!\!\!\!\!\!=0\!\!\quad\text{for $\beta=1,\dots,\dirnum$}. \label{abstract model nonlinear fokker--planck second line}
\end{align}
\end{subequations}
The meaning of the boundary conditions \eqref{abstract model nonlinear fokker--planck second line} is that the PDE \eqref{abstract model nonlinear fokker--planck first line} is satisfied in the weak sense when tested against any sufficiently smooth function $\psi$ satisfying $\nabla_u\psi(u)\cdot\textbf{n}_{\partial\R^{\dirnum}_+}\equiv 0$.

Under suitable assumptions, Theorems 2.4 and 2.5 in \cite{carrillo_clini_solem_the_mean_field_limit} ensure strong existence and uniqueness for the McKean--Vlasov equation \eqref{abstract 4 model McKean--Vlasov} and that $f(t,x,du)\coloneqq\law_{\R^{\dirnum}}(\Bar{u}(x,t))$ is the unique weak solution of the Fokker--Planck PDE \eqref{abstract model nonlinear fokker--planck}.

One of the key features of the model is that equation \eqref{abstract model nonlinear fokker--planck} has a unique solution independent of $\epsilon>0$.
That is to say, the law $f(t,x,du)=\law_{\R^{\dirnum}}(\Bar{u}^{\epsilon}(y,t))$ of a single McKean--Vlasov particle is independent of the noise correlation radius $\epsilon$.
This is because in \eqref{abstract 4 model McKean--Vlasov} the interaction between two neurons $\Bar{u}^{\epsilon}(x,t)$ and $\Bar{u}^{\epsilon}(y,t)$ happens only through their law $f$ as single random variables. 
In fact we will see that the dependence on the correlation radius is recovered as soon as we consider the joint law $f(x,y,t,,du,dv)=\law_{\R^{2\dirnum}}(\Bar{u}^{\epsilon}(x,t),\Bar{u}^{\epsilon}(y,t))$ of two or more particles (cf. Remark \ref{remark/dependence on epsilon}).

We also remark that the Fokker--Planck PDE \eqref{abstract model nonlinear fokker--planck}, in the case of diffusion term $\sigma\equiv \sigma_0$ constant, has been further analyzed in \cite{CHS2020,carrillo_roux_solem_noise_driven_bifurcation}.
Under particular assumptions on the drift term, it is shown that for every noise strength $\sigma_0>0$ there exists a unique stationary space homogeneous solution (from the modelling point of view this corresponds to the animal completely loosing track of its position in space).
This state is asymptotically stable for noise $\sigma_0$ big enough, but phase transition with multiple bifurcations, whose shape is also characterized locally, occurs as $\sigma_0$ decreases.

Finally, in \cite{carrillo_clini_solem_the_mean_field_limit} it is shown that the empirical measure $f_{MN}$ of the particle system \eqref{abstract model particle system} converges in a suitable sense (see Theorem \ref{theorem/ convergence of the single and joint empirical measure} below) to the weak solution $f$ of the Fokker--Planck equation \eqref{abstract model nonlinear fokker--planck}.
Informally, this yields the zeroth order expansion, as elements of $\pr(Q\times\R^{\dirnum})$,
\begin{equation}
    \label{formula/ zeroth order expansion}
    f_{MN}=f+O\left(\frac{1}{\sqrt{M}}+\frac{1}{N^{\nicefrac{\alpha}{d}}}\right)\quad\text{as $M,N\to\infty$}.
\end{equation}

The aim of the present work is to sharpen this picture by proving a central limit theorem for the fluctuations of the empirical measure around the deterministic limit $f$.
Precisely, we consider the rescaled fluctuations $\eta_t^{MN}:=\sqrt{M}(f_{MN}-f)$.
Proposition \ref{porposition/ real semimartingale rewriting of etaMN} below proves that, for suitable test functions $\psi$, the fluctuations admits the semimartingale expression
\begin{equation}
    \label{formula/ real semimartingale rewriting of etaMN / introduction}
    \langle\eta_t^{MN},\psi\rangle =  \langle\eta_0^{MN},\psi\rangle + \int_0^t  \langle\eta_r^{MN},\lin_r(f,f_{MN})[\psi]\rangle\,dr + M_t^{MN}(\psi),
\end{equation}
where 
\begin{equation}
    \label{formula/ real martingale term explicit expression / introduction}
    M_t^{MN}(\psi)=\frac{\sqrt{M}}{MN}\sum_{i=i}^N\sum_{k=1}^M\int_0^t\nabla_u\psi(x_i,u_{ik}(r))\sigma(x_i,r, u_{ik}(r),f_{MN})\, dW_{ik}(r)
\end{equation}
is a martingale term coming from the stochasticity of $f_{MN}$ and $\lin(f,f_{MN})$ is a suitable differential operator, essentially the linearization of the generator of the SDE \eqref{abstract model particle system} around $f$.
The aforementioned result already ensures that $f_{MN}\to f$ as $M,N\to\infty$.
A formal central limit theorem argument suggests also that $\eta_0^{MN}\to\eta_0$ for some Gaussian random variable in a suitable distribution space.
Finally, by analyzing the quadratic variation of the martingale $M_t^{MN}$ (cf. Proposition \ref{proposition/ pointwise convergence of the quad variation} below) we expect this to converge to some Gaussian process $G_t^{\epsilon}$ in a suitable distribution space.
In turn, from \eqref{formula/ real semimartingale rewriting of etaMN / introduction}, we conjecture that $\eta_t^{MN}$ converges to the solution of the Langevin SPDE (in integral form)
\begin{equation}
    \label{formula/ langevin spde / introduction}
    \eta_t^{\epsilon}=\eta_0+\int_0^t\lin_r(f,f)^*[\eta_r^{\epsilon}]\,dr + G_t^\epsilon,
\end{equation}
where the operator $\lin_r(f,f)$ is basically the linearization of the generator of the McKean--Vlasov SDE \eqref{abstract 4 model McKean--Vlasov} around $f$, and $\eta_0$ and $G_t^{\epsilon}$ are the above mentioned Gaussian terms.
The aim of this paper is to make this ansatz rigorous.
Here $\alpha\in(0,1]$ denotes the H\"older regularity in space of the coefficients $b$ and $\sigma$ in \eqref{abstract model particle system} (cf. Section \ref{section/Hypotheses and notations}).

\begin{theorem*}[Theorem \ref{theorem/ clt/ clt section} below]
For every $\epsilon>0$ fixed, along a scaling regime $M,N\to\infty$ such that $\sqrt{M}N^{-\nicefrac{\alpha}{d}}\to0$, the fluctuations $\eta_t^{MN}=\sqrt{M}(f_{MN}-f)$ of the empirical measure $f_{MN}$ of the particle system \eqref{abstract model particle system} around its deterministic limit, the solution $f$ of the Fokker--Planck equation \eqref{abstract model nonlinear fokker--planck}, satisfy
\begin{equation}
    \label{formula/ clt / introduction}
    \eta_t^{MN}\to\eta^{\epsilon}_t\quad\text{in law in $C\left([0,T];\mathcal{H}\right)$,}
\end{equation}
where $\eta^\epsilon_t$ is the unique solution of the Langevin SDPE \eqref{formula/ langevin spde / introduction} and $\mathcal{H}$ is a suitable Hilbert space of distributions on $Q\times \R^{\dirnum}_+$.
\end{theorem*}

The decay involving $M^{-\nicefrac{1}{2}}$, related to the mean field interaction along columns, and $N^{-\nicefrac{\alpha}{d}}$, related to the correlation between different columns, already appears in the mean field expansion \eqref{formula/ zeroth order expansion}.
The rescaling factor $\sqrt{M}$ in front of the fluctuations corresponds to the classical convergence rate for a mean field problem with $M$ particles subjected to independent sources of noise and is needed so that the fluctuations do not to vanish in the limit (cf. \cite{hitsuda_mitoma_tightness_problem_and_stochastic,Ferland_Fernique_Giroux_compactness_of_the_fluctuations,Bezandry_Fernique_Giroux_A_functional_central_limit_theorem, fernandez_meleard_a_hilbertian_approach}).
However, as discussed in Remark \ref{remark/n mean field problem}, we actually have $N$ correlated `$M$-particle mean field clusters' interacting among themselves with strength $N^{-\nicefrac{\alpha}{d}}$.
The scaling regime $\sqrt{M}N^{-\nicefrac{\alpha}{d}}\to0$ serves to prevent these clusters from resonating and resulting in yet unbounded fluctuations (cf. Remark \ref{remark/ rescaling factor sqrt M}).

In the same fashion as \eqref{formula/ zeroth order expansion}, the theorem above yields the first order expansion
\begin{equation}
    \label{formula/ first order expansion}
    f_{MN}=f+\frac{1}{\sqrt{M}}\eta_t^\epsilon+o\left(\frac{1}{\sqrt{M}}+\frac{1}{N^{\nicefrac{\alpha}{d}}}\right)\quad\text{as $M,N\to\infty$ and $\sqrt{M}N^{-\nicefrac{\alpha}{d}}\to0$.}
\end{equation}
This is particularly relevant from the simulation point of view: in order to understand the behavior of a large number of grid cells obeying the model \eqref{abstract model particle system}, simulating the continuous model $f+M^{-\nicefrac{1}{2}}\eta_t^\epsilon$ is numerically much cheaper than simulating $MN$ interacting SDEs.
On the other hand, our main interest is the effect of noise on grid cells.
The zero order term $f$ alone completely loses track of the noise present at the particle level, which is instead retained by the first order expansion \eqref{formula/ first order expansion} thus furnishing a sharper description.

\begin{remark}
\label{remark/ joint scaling regime in m n epsilon}
We point out that we could also let $\epsilon=\epsilon(M,N)$ vary with the number of neurons and consider a joint scaling regime $M,N\to\infty$ and $\epsilon(M,N)\to\Bar{\epsilon}$.
As long as $\Bar{\epsilon}>0$, all the results of the paper hold with virtually identical proofs.
Some results and estimates also extend to the case $\epsilon(M,N)\to0$, possibly after imposing some constraint on the scaling regime between $\epsilon(M,N)\to0$ and $N\to\infty$.
However, the point is that in the case $\epsilon(M,N)\to0$, the limiting objects for the fluctuations and, even before, the supposed mean field limit \eqref{abstract 4 model McKean--Vlasov} for the SDE system \eqref{abstract model particle system} is mathematically ill-posed.
For $\epsilon=0$ the noise $W^{\epsilon}(x,t)$ featuring in the McKean--Vlasov equation \eqref{abstract 4 model McKean--Vlasov} makes sense only as a distribution (see e.g. \cite[Chapter 4]{daprato_zabczyk_1992}) and in turn the equation, even in the simpler case $\sigma(x,t,u,f)\equiv \sigma_0$ constant, involves applying the drift nonlinearity $b(\dots)$ to a distribution.
It is not clear whether it is even possible to give a meaning to equation \eqref{abstract 4 model McKean--Vlasov} for $\epsilon=0$ (see also \cite[Remark 1.2]{carrillo_clini_solem_the_mean_field_limit}) and in turn to the joint law $f^{2,\epsilon}$, introduced in \eqref{formula/joint law of mckean--vlasov particles} below, featuring in the definition \eqref{formula/ quadratic variation of real gaussian limit} of the quadratic variation of the fluctuation limit $\eta_t^\epsilon$.
\end{remark}


\bigskip

\noindent
\emph{\textbf{Overview of the methods.}}

The methods of this paper are deeply based on the approach to fluctuations put forward in Fernandez and Méléard \cite{fernandez_meleard_a_hilbertian_approach}, Bezandry, Ferland, Fernique and Giroux \cite{Bezandry_Fernique_Giroux_A_functional_central_limit_theorem,Ferland_Fernique_Giroux_compactness_of_the_fluctuations} and Hitsuda and Mitoma \cite{hitsuda_mitoma_tightness_problem_and_stochastic}.
Our situation presents additional difficulties: unboundedness of the coefficients, boundary conditions imposing positivity of the activity level of the neurons, and interaction between different mean-field families of coupled neurons across the cortex. 

The fluctuations $\eta_t^{MN}=\sqrt{M}(f_{MN}-f)$ are random signed measures in the path space over $Q\times\R^{\dirnum}$, whose space part is independent both of the time and the stochasticity.
The first problem is to find a distribution space in which both $\eta_t^{MN}$ and its limit belong.
In Theorem \ref{theorem/ clt/ clt section} we prove that, along the scaling regime $\sqrt{M}N^{-\nicefrac{\alpha}{d}}\to0$, the fluctuations converge in law in $C\big([0,T];H_x^{-\eplus}\h^{-\ktwo,\ttwo}\big)$.
The spaces $H_x^{e}\h^{\kappa,\theta}$, introduced in Section \ref{section/Hypotheses and notations}, are vector valued Sobolev spaces taking values in other weighted Sobolev spaces with suitable boundary conditions needed to account for the reflecting boundary conditions at the SDE level, which reduce the class of admissible test functions, and the space interaction among neurons.
The particular exponents $e$, $\kappa$ and $\theta$, given in \eqref{formula/ values of kappa theta e}, depend on the dimension of the spaces $Q\subseteq\R^d$ and $\R^{\dirnum}$ to furnish suitable Sobolev embeddings.
It will appear from the proofs that the space $H_x^{-\eplus}\h^{-\ktwo,\ttwo}$ is optimal and provides the minimal Hilbert space in which to embed the fluctuations (cf. Remark \ref{remark/ on the need for hilbert setting vs holderianity only}).

Our strategy revolves around the semimartingale equation \eqref{formula/ real semimartingale rewriting of etaMN / introduction} for $\eta_t^{MN}$.
Estimates on the norm of $M_t^{MN}$ are readily obtained with martingale arguments (cf. Proposition \ref{proposition/ H valued martingale MtMN}).
Here the main difficulty is that the operator $\lin(f,f_{MN})$ featuring in \eqref{formula/ real semimartingale rewriting of etaMN / introduction} is a second order differential operator and thus not bounded in any particular space of the form $H_x^{e}\h^{\kappa,\theta}$, since it reduces the regularity of the test functions.
Hence we cannot use equation \eqref{formula/ real semimartingale rewriting of etaMN / introduction} and Gr\"onwall's Lemma to estimate the norm of $\eta_t^{MN}$.
To circumvent this problem, we work between two nested spaces $H_x^{\eplus}\h^{\ktwo,\ttwo}\subseteq H_x^{\e}\h^{\kone,\tone}$ where the operator $\lin(f,f_{MN})$ results bounded (cf. Lemma \ref{proposition/ norm bounds on the linearized operator}) and we obtain the needed estimates on $\eta_t^{MN}$ (cf. Lemma \ref{proposition/ first bounds on etaMN} and Proposition \ref{proposition/ H valued semimartingale etaMN}).
It is the combination of the estimates for $M_t^{MN}$ and $\eta_t^{MN}$ that imposes the particular scaling regime $\sqrt{M}N^{-\nicefrac{\alpha}{d}}\to0$ (cf. Remark \ref{remark/ rescaling factor sqrt M}).

These estimates are readily translated into tightness of the processes $M_t^{MN}$ and $\eta_t^{MN}$ thanks to the compactness of the embedding $H_x^{\eplus}\h^{\ktwo,\ttwo}\subseteq H_x^{\e}\h^{\kone,\tone}$ and martingale methods (cf. Proposition \ref{proposition/ tightness of the processes}).
A fine analysis based on Sznitman's coupling method \cite{Sznitman1984} shows that the quadratic variation of the martingale term $M_t^{MN}$ is converging to a deterministic limit (cf. Proposition \ref{proposition/ pointwise convergence of the quad variation}) and this uniquely identifies a Gaussian process $G_t^{\epsilon}$ such that $M_t^{MN}\to G_t^{\epsilon}$ (cf. Proposition \ref{proposition/ clt for the martingale}).
An adaptation of the standard Lévy central limit theorem shows the convergence of the initial data $\eta_0^{MN}$ (cf. Lemma \ref{proposition/ clt for the initial data}).
Theorem \ref{theorem/ clt/ clt section} follows from passing to the limit in equation \eqref{formula/ real semimartingale rewriting of etaMN / introduction} and the well-posedness of the Langevin SPDE \eqref{formula/ langevin spde / introduction}.


\bigskip

\noindent
\emph{\textbf{Structure of the work.}}

The work is structured as follows. In the next section we first present our hypotheses and notation and deduce some immediate consequences for later convenience. 
Then we lay out the precise setting for our results and we collect some auxiliary results needed in our arguments.
In Section \ref{section/the central limit theorem} we establish our central limit theorem following the strategy discussed above.
In Section \ref{section/proof of the results} we collect the proofs of the results.



\section{Preliminaries and auxiliary results}

\subsection{Hypotheses and notation}
\label{section/Hypotheses and notations}

In this section we introduce our hypotheses and notations.
We fist take care of the initial data and the noise terms in the McKean--Vlasov equation \eqref{abstract 4 model McKean--Vlasov}.
The initial data and the noise in the particle system \eqref{abstract model particle system} will then be given by the corresponding quantities at the locations $x_i\in Q$ (cf. Setting \ref{setting}).
As regards the initial data, we will consider random families $(u(x,0))_{x\in Q}$ such that 
\begin{equation}\label{formula/ assumption on initial data}
    u(\cdot,0)\in C^{\alpha}\big(Q;L^4(\Omega)\big)\cap L^{\infty}\big(Q;L^{4\tone}(\Omega)\big)\,\, \text{for some }\alpha\in(0,1]\text{ and } \tone>1+\frac{\dirnum}{2}+\max\left\{2,\frac{\dirnum}{2}\right\}.
\end{equation}
The H\"older continuity in space is needed to control the interaction among columns at different locations and is motivated by the modelling, as we expect the randomness to be somewhat continuous along the cortex.
The moments will serve to obtain bounds for the fluctuations in the dual space of suitably regular function spaces (cf. \eqref{formula/ norm of evaluation operators}).

As regards the noise, we consider a Gaussian random field $W^{\epsilon}:\Omega\times\R^d\times\R_+\to\R^{\dirnum}$ with independent components $\beta=1\dots\dirnum$, zero mean and covariance 
\begin{equation}\label{formula/ space colored time white noise}
    \E\left[W^{\epsilon,\beta}(x,t)W^{\epsilon,\beta}(y,s)\right]\!=\!(t\wedge s)\, C_{\rho}\,\epsilon^d\int_{\R^d}\rho_\epsilon(z-x)\rho_{\epsilon}(z-y)\,dz,\quad \text{for } C_{\rho}=\left(\int_{\R^d}\rho(z)^2\,dz\right)^{-1}\!\!,
\end{equation}
where $\rho:\R^d\to[0,1]$ is a radial mollifier supported in the unitary ball and $\rho_\epsilon$ is its the $\epsilon$-rescaled version.
Such a process can be obtained for example from convolution and rescaling of independent copies $W^\beta$ of an actual distribution valued space-time white noise (see e.g. \cite[Chapter 4]{daprato_zabczyk_1992}), by setting
\begin{equation}\label{formula/epsilon convoluted and rescaled noise}
    W^{\epsilon,\beta}(x,t)\coloneqq C_{\rho}^{\frac{1}{2}}\,\epsilon^{\frac{d}{2}}\,\langle W_t,\rho_{\epsilon}(\cdot-x)\rangle.
\end{equation}
As mentioned in Remark \ref{remark/ joint scaling regime in m n epsilon}, the $\epsilon$-correlation in space is needed to avoid mathematical issues and in fact, as for the initial data, is strongly motivated by the modelling.

For future reference, we highlight some of the properties of $W^\epsilon$.
First, from \eqref{formula/ space colored time white noise} we have that $\E[W^{\epsilon,\beta}(x,t)W^{\epsilon,\gamma}(x,s)]=\delta_{\beta\gamma}\,\,t\wedge s$. 
Thus, for fixed $x$, the process $t\mapsto W^{\epsilon}(x,t)$ is a $\dirnum$-dimensional Brownian motion.
Similarly, from $\sup(\rho)\subseteq B(0,1)$ it follows that
\begin{equation}
    \E\left[W^{\epsilon}(x,t)W^{\epsilon}(y,s)\right]=0 \,\,\text{ if }\,\, |x-y|>2\epsilon.
\end{equation}
Hence the processes $W^{\epsilon}(x,t)$ and $W^{\epsilon}(y,t)$ are independent for $|x-y|>2\epsilon$.
Furthermore, using \eqref{formula/ space colored time white noise}, we compute
\begin{align}
    \E\left[\left|W^{\epsilon}(x,t)-W^{\epsilon}(y,s)\right|^2\right]
    &\leq C\,
    \E\left[\left|W^{\epsilon}(x,t)-W^{\epsilon}(x,s)\right|^2+\left|W^{\epsilon}(x,s)-W^{\epsilon}(y,s)\right|^2\right]
    \\
    &\leq C\left(
    |t-s|+s\, C_\rho\, \epsilon^d\int_{\R^d}\left(\rho_\epsilon(z-x)-\rho_{\epsilon}(z-y)\right)^2\,dz\right)
    \\
    &\leq C
    \left(|t-s|+\frac{|x-y|^2}{\epsilon^2}\right),
\end{align}
for a constant $C=C(\rho)$.
Similar estimates hold for any higher moment $p\geq 2$ and the Kolmogorov continuity theorem ensures the existence of a suitable modification of $W^\epsilon$ with continuous trajectories in both $x$ and $t$.
In particular, we have that $W^{\epsilon}$ is jointly measurable in $(x,t)\in\R^d\times\R^+$ and in the sample path $\omega\in\Omega$.
Finally, for any $x,y\in\R^d$, a direct computation shows that the quadratic variation of the martingale $W^{\epsilon,\beta}(x,t)-W^{\epsilon,\beta}(y,t)$ satisfies
\begin{gather}\label{formula/epsilon convoluted noise quadratic variation}
\begin{aligned}
    \big\langle W^{\epsilon,\beta}(x,t-W^{\epsilon,\beta}(y,t)\big\rangle=t\,C_\rho\, \epsilon^d\int_{\R^d}\left(\rho_\epsilon(z-x)-\rho_{\epsilon}(z-y)\right)^2\,dz
    \leq t\,
    C\,\frac{|x-y|^2}{\epsilon^2},
\end{aligned}
\end{gather}
for a constant $C=C(\rho)$.

Next we present the assumptions on the coefficients of the equations.
We remark that the specific shape and properties of these generalized coefficients are deduced from those of the concrete model analyzed in \cite{carrillo_clini_solem_the_mean_field_limit}.
Namely, the coefficients $b,\sigma:Q\times\R^+\times\R^{\dirnum}\times\pr(Q\times\R^{\dirnum})\to\R^{\dirnum}$ take the form
\begin{align}\label{formula/ shape of coefficients}
\begin{split}   &b^\beta(x,r,u,\mu)=b_0^\beta(x,r,u)+\phi^\beta\left(\int_{Q\times\R^{\dirnum}}b_1^\beta(x,y,r,u,v)\,\mu(dy,dv)\right),
\\
&\sigma^\beta(x,r,u,\mu)=\sigma_0^\beta(x,r,u)+\phi^\beta\left(\int_{Q\times\R^{\dirnum}}\sigma_1^\beta(x,y,r,u,v)\,\mu(dy,dv)\right),
\end{split}
\end{align}
for suitable functions $b_0,\sigma_0:Q\times\R^+\times\R^{\dirnum}\to\R^{\dirnum}$ and $b_1,\sigma_1:Q^2\times\R^+\times\R^{2\dirnum}\to\R^{\dirnum}$, and $\phi:\R\to\R^{\dirnum}$.
\begin{remark}\label{remark/ same coefficients for every direction}
In order simplify the exposition and the computations, we shall assume the coefficients are the same for every direction.
That is to say, for every $\beta=1,\dots,\dirnum$,
\begin{align}
\begin{split}   b^\beta(x,r,u,\mu)=b_0(x,r,u)+\phi\left(\int_{Q\times\R^{\dirnum}}\!\!\!\!\!b_1(x,y,r,u,v)\,\mu(dy,dv)\right),
\end{split}
\end{align}
and similarly for the coefficients $\sigma$, for suitable functions $b_0$, $b_1$, $\sigma_0$, $\sigma_1$ and $\phi$ satisfying the assumptions below.
We stress that all the results hold identically with identical proofs for general direction dependent coefficients $b^\beta$, $\sigma^\beta$ and $\phi_\beta$ satisfying the assumptions below for each single $\beta$.
Furthermore, to ease the notation, we will also use the shorthand
\begin{equation}
\label{formula/ notation for b1(x,t,u,mu)}
b_1(x,t,u,\mu):=\int_{Q\times\R^{\dirnum}}\!\!\!\!b_1(x,y,t,u,v)\, \mu(dy,dv),\quad \sigma_1(x,t,u,\mu):=\int_{Q\times\R^{\dirnum}}\!\!\!\!\sigma_1(x,y,t,u,v)\, \mu(dy,dv).
\end{equation}
\end{remark}

We assume the nonlinearity $\phi$ is smooth with bounded derivatives.
Precisely we assume that 
\begin{equation}\label{formula/ assumption on phi}
    \frac{d^j\!\phi}{dz^j}(z)\in L^{\infty}(\R) \quad\forall\,j\geq1.
\end{equation}
In particular $\phi$ is globally Lipschitz and with sublinear growth.
We assume that $b_0$ and $\sigma_0$ are smooth in the $x$ and $u$ variables, uniformly $\alpha$-H\"older continuous in the $x$-variable (for $\alpha\in(0,1]$ given in \eqref{formula/ assumption on initial data}) and with $u$-derivatives bounded uniformly in $x$ and $t$.
That is to say, for every $h,k\in\N$, for suitable constants $C$,
\begin{align}\label{formula/ assumption on b0 sigma0 first}
    \begin{split}
        |D_x^hD_u^{k+1}b_0(x,t,u)|&\leq C, 
        \\ 
        \left|D_x^hb_0(x,t,u)-D_x^hb_0(x',t,u)\right|&\leq C |x-x'|^\alpha,
    \end{split}
    \,\,\, \forall x,x',t,u\in Q^2\times\R_+\times\R^{\dirnum},
\end{align}
and analogously for $\sigma_0$.
In particular this implies Lipschitzianity and sublinear growth in $u$. That is, for every $h,k\in\N$, for suitable constants $C$,
\begin{align}\label{formula/ assumption on b0 sigma0 second}
    \begin{split}
        \left|D_x^hb_0(x,t,u)\right|\leq C(1+|u|),
        \\ 
        \left|D_x^hD_u^kb_0(x,t,u)-D_x^hD_u^kb_0(x,t,u')\right|&\leq C |u-u'|,     \end{split}
    \,\,\, \forall x,t,u,u'\in Q\times\R_+\times\R^{2\dirnum},
\end{align}
and analogously for $\sigma_0$.
Similarly, we assume that $b_1$ and $\sigma_1$ are smooth in the $x,y$ and $u,v$ variables, uniformly $\alpha$-H\"older in the $x,y$-variable and with $u,v$-derivatives bounded uniformly in $x,y$ and $t$.
That is, for every $h,k\in\N$, for suitable constants $C$, for every $x,x',y,y',t,u,v\in Q^4\times\R_+\times\R^{2\dirnum}$, 
\begin{align}\label{formula/ assumption on b1 sigma1 first}
    \begin{split}
        |D_{x,y}^hD_{u,v}^{k+1}b_1(x,y,t,u,v)|&\leq C, 
        \\ 
        \left|D_{x,y}^hb_1(x,y,t,u,v)-D_{x,y}^hb_0(x',y',t,u,v)\right|&\leq C |x-x'|^\alpha+|y-y'|^\alpha,
    \end{split}
\end{align}
and analogously for $\sigma_1$.
In particular this implies Lipschitzianity and sublinear growth in $u$ and $v$.
That is, for every $h,k\in\N$, for suitable constants $C$, for every $x,y,t,u,u',v,v'\in Q^2\times\R_+\times\R^{4\dirnum}$, 
\begin{align}\label{formula/ assumption on b1 sigma1 second}
    \begin{split}
        \left|D_{x,y}^hb_1(x,y,t,u,v)\right|\leq C(1+|u|+|v|),
        \\ 
        \left|D_{x,y}^hD_{u,v}^kb_1(x,y,t,u,v)-D_{x,y}^hD_{u,v}^kb_1(x,y,t,u',v')\right|&\leq C (|u-u'|+|v-v'|),    \end{split}
\end{align}
and analogously for $\sigma_1$.

The shape \eqref{formula/ shape of coefficients} of the coefficients and assumptions \eqref{formula/ assumption on phi}-\eqref{formula/ assumption on b1 sigma1 second} immediately translate into the following properties of $b$ and $\sigma$.

\begin{lemma}\label{proposition/ properties of b and sigma}
The coefficients $b,\sigma:Q\times\R^+\times\R^{\dirnum}\times\pr(Q\times\R^{\dirnum})\to\R$ are smooth in the $x$ and $u$ variables.
Furthermore we have, for every $h,k\in\N$, for suitable constants $C\geq0$, for every $x,x',t,u,u',\mu,\mu'\in Q^2\times\R_+\times\R^{2\dirnum}\times\pr(Q\times\R^{\dirnum})^2$, 
\begin{align}\label{formula/ properties of b and sigma}
    \begin{split}
        &|D_x^hD_u^{k+1}b(x,t,u,\mu)|\leq C,
        \\
        &|D_x^hb(x,t,u,\mu)|\leq C\left(1+|u|+\int_{Q\times\R^{\dirnum}}\!\!\!|v|\,\mu(dy,dv)\right), 
        \\
        &|D_x^hD_u^kb(x,t,u,\mu)-D_x^hD_u^kb(x',t,u',\mu')|
        \\
        &\qquad\leq C\left(|x-x'|^\alpha+|u-u'|+\w_1(Q\times\R^{\dirnum})(\mu,\mu')^\alpha+\w_1(Q\times\R^{\dirnum})(\mu,\mu')\right),
    \end{split}
\end{align}
and analogously for $\sigma$, where $\w$ denotes the Wasserstein distance (see formula \eqref{formula/ wasserstein distance} below).
\end{lemma}

We now introduce the function spaces needed to state our results.
First of all, given a Banach space $X$, we consider the space $\pr_m(X)$ of probability measures on $X$ with finite $m^{th}$ moment, which is a Polish space when endowed with the $m^{th}$  order Wasserstein distance
\begin{equation}\label{formula/ wasserstein distance}
    \w_m(X)(P,Q):=\inf\left\{\int_{X^2}|x-y|^m\,d\pi\mid\,\text{$\pi$ pairing between $P$ and $Q$}\right\}^{\nicefrac{1}{m}}.
\end{equation}

Next we introduce the following weighted Sobolev spaces with no flux boundary conditions, for any integer $k\geq 2$ and any $\theta\geq0$,
\begin{equation}
    \h^{k,\theta}:=\left\{\psi\in W_{loc}^{j,2}(\R^{\dirnum}_+)\mid \nabla\psi{(u)}\cdot\textbf{n}_{\partial\R^{\dirnum}_+}\equiv0,\quad \sum_{j=0}^{k}\int_{\R^{\dirnum}_+}\frac{|D_u^j\psi(u)|^2}{1+|u|^{2\theta}}\,du<\infty\right\},
\end{equation}
which are Hilbert spaces when endowed with the norm $\|\psi\|_{\h^{k,\theta}}=\left(\sum_{j=0}^{k}\int_{\R^{\dirnum}_+}\frac{|D_u^j\psi(u)|^2}{1+|u|^{2\theta}}\,du\right)^{\nicefrac{1}{2}}$.
Similarly, we will also consider the Banach spaces
\begin{equation}
    \mathscr{W}^{k,\infty,\theta}:=\left\{\psi\in W_{loc}^{k,\infty}(\R^{\dirnum}_+)\mid \nabla\psi(u)\cdot\textbf{n}_{\partial\R^{\dirnum}_+}\equiv0,\quad \sup_{u\in\R^{\dirnum}_+}\sum_{j=0}^{k}\frac{|D_u^j\psi(u)|}{1+|u|^{\theta}}<\infty\right\},
\end{equation}
endowed with the analogous norm $\|\psi\|_{\mathscr{W}^{k,\infty,\theta}}=\sup_{u\in\R^{\dirnum}_+}\sum_{j=0}^{k}\frac{|D_u^j\psi(u)|}{1+|u|^{\theta}}$.

\begin{remark}\label{remark/ no flux boundary conditions}
These spaces are particularly adapted to our needs because of the no-flux condition $\nabla\psi(u)\cdot\textbf{n}_{\partial\R^{\dirnum}_+}\equiv0$.
This ensures that $\psi$ is an admissible test function for the Fokker-Planck equation \eqref{abstract model nonlinear fokker--planck} according to our definition of weak solution.
Furthermore it ensures that, when applying It\^o formula to $\psi(u_{ik}(t))$ for $u_{ik}$ solving \eqref{abstract model particle system} or to $\psi(\Bar{u}^\epsilon(x,t))$ for $\Bar{u}^\epsilon(x,t)$ solving \eqref{abstract 4 model McKean--Vlasov} respectively, the resulting reflection term 
\begin{equation}
    \int_0^t\nabla\psi(u_{ik}(r))\,d\ell_{ik}(r)= -\int_0^t\nabla\psi(u_{ik}(r))\cdot \textbf{n}_{\partial\R^{\dirnum}_+}\,1_{\{u_{ik}\in\partial\R^{\dirnum}_+\}}\,d|\ell_{ik}|(r)=0 ,
\end{equation}
or $\int_0^t\nabla\psi(\Bar{u}^\epsilon(x,r)) d\ell(x,r)=0$ respectively, is identically zero.  
\end{remark}

Classical embedding results extend to these particular Sobolev spaces (see e.g. \cite{adams_sobolev_spaces}).
In particular, we will make use of the following \emph{continuous} inclusions
\begin{equation}\label{formula/ morrey embedding of normal weighted sobolev spaces}
    \h^{k+j,\theta}\hookrightarrow\mathscr{W}^{j,\infty,\theta}\quad \forall\, k>\frac{\dirnum}{2}\,\,\,\forall\,j\geq0\,\,\,\forall\,\theta\geq0,
\end{equation}
and the following \emph{compact} inclusions
\begin{equation}\label{formula/ compact embedding of normal weighted sobolev spaces}
    \h^{k+j,\theta}\doublehookrightarrow\h^{j,\theta+\beta}\quad \forall\, k>\nicefrac{\dirnum}{2}\,\,\,\forall\, \beta>\nicefrac{\dirnum}{2}\,\,\,\forall\,j\geq0\,\,\,\forall\,\theta\geq0.
\end{equation}

Given a Hilbert space $V$ we will also consider the vector valued H\"older and Sobolev spaces, for some $e\geq0$,
\begin{equation}
    C_x^\alpha V:=C^{\alpha}(Q;V),\quad H_x^eV:=W^{e,2}(Q;V).
\end{equation}
In particular, $H_x^eV$ is a Hilbert spaces with the usual Sobolev norm (cf. Remark \ref{remark/ on the need for hilbert setting vs holderianity only} below).
We will need the following vector valued versions of Sobolev embeddings (see e.g. \cite{amann_compact_embedding_vector_valued_sobolev}).
For any Hilbert space $V$ we have the continuous embedding
\begin{equation}\label{formula/ continuous embedding of vector valued sobolev spaces}
    H_x^eV\hookrightarrow C_x^\alpha V\quad\forall\, e>\alpha+\nicefrac{d}{2}.
\end{equation}
For any compact embedding $V\doublehookrightarrow V_0$ of Hilbert spaces, we have the compact embedding
\begin{equation} \label{formula/ compact embedding of vector valued sobolev spaces}
    H_x^{e^+}V\doublehookrightarrow H_x^e V_0\quad\forall\,e\geq0,
\end{equation}
where we denoted $e^+:=e+\gamma$ for any arbitrary $\gamma>0$.

In particular, combining \eqref{formula/ morrey embedding of normal weighted sobolev spaces}-\eqref{formula/ compact embedding of vector valued sobolev spaces} we obtain the following embeddings, which we collect here for later convenience,
\begin{align}
    H_x^e\h^{k+j,\theta}&\hookrightarrow C_x^{\alpha}\mathscr{W}^{j,\infty,\theta}\text{ continuously }\quad \forall e>\frac{d}{2}+\alpha\,\,\,\forall\,k>\frac{\dirnum}{2}\,\,\,\forall j,\theta\geq0, \label{formula/ continuous embedding of working spaces into bounded functions}
    \\
    H_x^{e^+}\h^{k+j,\theta}&\doublehookrightarrow H_x^e\h^{j,\theta+\beta} \text{ compactly }\quad \forall e\geq 0\,\,\,\forall\,k,\beta>\frac{\dirnum}{2}\,\,\,\forall j,\theta \geq0. \label{formula/ compact embedding between working spaces}
\end{align}

\begin{remark}\label{remark/ on the need for hilbert setting vs holderianity only}
The spaces $H_x^e\h^{\kappa,\theta}$ are introduced for the sole reason of working with subspaces of $C^\alpha(Q;\h^{\kappa,\theta})$ possessing a Hilbert space structure so as to simplify many arguments in the following.
In fact, results completely analogous to those presented could be obtained by working directly in the Banach spaces $C^\alpha(Q;\h^{\kappa,\theta})$.
However, having no Hilbert space structure at disposal, the methods of this paper need to be carefully adapted by introducing several technicalities which, in the author's opinion, would not improve drastically the results stated.
To start with, even to state the results, we would have to introduce the theory of stochastic integration in Banach spaces, put forward only in recent years (see e.g. \cite{neerven_stochastic_integration_banach_spaces}).
\end{remark}

Finally, we introduce the following `evaluation operators' and we estimate their norm for later use.
Namely, for any $j\in\N$ and any $(x,u),(y,v)\in Q\times\R^{\dirnum}_+$, for $a=\left(\substack{\dirnum+j-1\\j}\right)$, we define
\begin{align}\label{formula/ evaluation operators}
    &V_{(x,u)}^j:C_x^\alpha \mathscr{W}^{j,\infty,\theta}\to\R^a\,\big|\,\psi\mapsto D_u^j\psi(x,u),
    \\
    &V^{j,\text{dif}}_{(x,u),(y,v)}:C_x^\alpha \mathscr{W}^{j+1,\infty,\theta}\to\R^a\,\big|\,\psi\mapsto \left(D_u^j\psi(x,u)-D_u^j\psi(y,v)\right).
\end{align}
Using the definition of the spaces $\mathscr{W}^{j,\infty,\theta}$ we compute
\begin{align}
    |V_{(x,u)}^j\psi|=|D_u^j\psi(x,u)|\leq
    \sup_{y,v}\frac{|D^j_u\psi(y,v)|}{1+|v|^\theta}\,(1+|u|^\theta)\lesssim\|\psi\|_{C_x^\alpha\mathscr{W}^{j,\infty,\theta}}\,(1+|u|^\theta),
\end{align}
and similarly
\begin{align}
    |V_{(x,u),(y,v)}^{j,\text{dif}}\psi|
    &\leq
    |D_u^j\psi(x,u)-D_u^j\psi(y,u)|+|D_u^j\psi(y,u)-D_u^j\psi(y,v)|
    \\
    &\lesssim
    \sup_{w}\frac{|D_u^j\psi(x,w)-D_u^j\psi(y,w)|}{1+|w|^\theta}\,(1+|u|^\theta)
    +
    \sup_{w\in\overrightarrow{uv}} \frac{|D_u^{j+1}\psi(y,w)|}{1+|w|^\theta}\,(1+|u|^\theta+|v|^\theta)\,|u-v|
    \\
    &\lesssim \|\psi\|_{C_x^\alpha\mathscr{W}^{j+1,\infty,\theta}}|x-y|^\alpha(1+|u|^\theta)
    +
    \|\psi\|_{C_x^\alpha\mathscr{W}^{j+1,\infty,\theta}}(1+|u|^\theta+|v|^\theta)|u-v|.
\end{align}
That is, we have bounds for $\|V_{(x,u)}^j\|_{(C_x^\alpha\mathscr{W}^{j,\infty,\theta})^*}$ and $\|V_{(x,u),(y,v)}^{j,\text{dif}}\|_{(C_x^\alpha\mathscr{W}^{j+1,\infty,\theta})^*}$.
Combining this with the embeddings \eqref{formula/ continuous embedding of working spaces into bounded functions} immediately yields the estimates, for suitable constants $C$,
{\small
\begin{align}\label{formula/ norm of evaluation operators}
\begin{split}
     \|V_{(x,u)}^j\|_{H_x^{-e}\h^{-(k+j),\theta}}&\leq C\,\left(1+|u|\right)^\theta,
     \\
     \|V_{(x,u),(y,v)}^{j,\text{dif}}\|_{H_x^{-e}\h^{-(k+j+1),\theta}}&\leq C\left(1+|u|+|v|\right)^\theta \big(|x-y|^\alpha+|u-v|\big),
\end{split}
\,\,\forall e>\frac{d}{2}+\alpha\,\,\forall k>\frac{\dirnum}{2}\,\,\forall j,\theta\geq0.
\end{align}
}


\subsection{Setting and auxiliary results} 
\label{section/setting and auxiliary results}

In this section we first lay out the precise setting for our results and then we collect some auxiliary results for later convenience.

\begin{notation*}
To ease the formulas for the function spaces involved, we will use the shorthands:
{\small
\begin{align}\label{formula/ values of kappa theta e}
\begin{aligned}
    &\kone:=\bigg(1+\frac{\dirnum}{2}\bigg)^+,
    &\ktwo:=\kone+\max\left\{2,\left(\frac{\dirnum}{2}\right)^+\right\},
    &\qquad\qquad\e=\alpha+\frac{d}{2},
    \\
    &\tone:=\ttwo+\max\left\{2,\left(\frac{\dirnum}{2}\right)^+\right\},
    &\ttwo:=\bigg(1+\frac{\dirnum}{2}\bigg)^+,\qquad\quad\qquad
    &
\end{aligned}
\end{align}
}
where for any $\beta\in\R$ we denote $\beta^+:=\beta+\gamma$ for any arbitrary $\gamma>0$.
\end{notation*}

This choice of Sobolev exponents and weights will give us the needed embeddings of type \eqref{formula/ continuous embedding of working spaces into bounded functions} and \eqref{formula/ compact embedding between working spaces}, and the needed bounds of type \eqref{formula/ norm of evaluation operators}.

\begin{setting}\label{setting}
We fix $\epsilon>0$ and a filtered probability space $(\Omega,(\mathcal{F}_t)_{t\geq0},\p)$, and we assume it supports all the random variables listed below.
First, for each $k\in\N$, let $\{W_k^\epsilon(x,t)\}_{k\in\N}$ be independent adapted $\dirnum$-dimensional Gaussian noise terms as in formula \eqref{formula/ space colored time white noise}.
Then, for $h\in\N$, we assume i.i.d. families of $\mathcal{F}_0$-measurable random initial conditions $u_h(x,0)\in C^{\alpha}\big(Q;L^4(\Omega)\big)\cap L^{\infty}\big(Q;L^{4\tone}(\Omega)\big)$ on the sheet $Q$. 
In particular, they are independent of the noise terms $\{W_k^\epsilon(x,t)\}_{k\in\N}$.
Next, for each $N\in\N$ we take points $x_1,\dots,x_N\in Q$ in the center of the squares of an equispaced grid on $Q=[0,1]^d$ with side length $N^{-\frac{1}{d}}$.
We denote by $Q^N_i$ the square with center $x_i$, and we notice that $\meas(Q_i^N)=\frac{1}{N}$ and $\text{diam}(Q_i^N)=\sqrt{d}N^{-\nicefrac{1}{d}}$.
Finally, for every $N,M\in\N$, we introduce the particles for the coupling method.
For $i=1,\dots,N$ and $k=1,\dots,M$, let $u_{ik}^{\epsilon}(t)$ be the solution of the particle system \eqref{abstract model particle system} with initial data $u_{ik}(0)\coloneqq u_k(x_i,0)$ and Brownian motions $W_{ik}^{\epsilon}(t)\coloneqq W^{\epsilon}_k(x_i,t)$.
For each $k\in\N$, let $\bar{u}_k^{\epsilon}(x,t)$ be the solution of the McKean--Vlasov equation \eqref{abstract 4 model McKean--Vlasov}  with initial data $u_k(x,0)$ and noise $W_k^{\epsilon}(x,t)$, and for $i=1,\dots,N$ define $\bar{u}_{ik}^{\epsilon}(x,t)\coloneqq \bar{u}^{\epsilon}_k(x_i,t)$.
\end{setting}

In this setting, Theorems 2.3, 2.4 and 2.5 in \cite{carrillo_clini_solem_the_mean_field_limit} ensure the systems \eqref{abstract model particle system} and \eqref{abstract 4 model McKean--Vlasov} are well-posed and that $f(x,t,du)=\law_{\R^{\dirnum}}(\Bar{u}_k^\epsilon(x,t))$ is the unique weak solution of the Fokker--Planck PDE \eqref{abstract model nonlinear fokker--planck} with initial data $f_0(x,du)=\law_{\R^{\dirnum}}(u_k(x,0))$.

We now collect some auxiliary results either new or taken from \cite{carrillo_clini_solem_the_mean_field_limit} we will need in the following.
First of all, owing the i.i.d. properties of the initial data and the noise terms, we have the following.
\begin{lemma}
The collections of particles $\{u^{\epsilon}_{ik}(t)\}_{i=1,\dots,N}$ are exchangeable families for $k=1,\dots,M$.
Moreover, the collections of particles $\{\bar{u}_{ik}^{\epsilon}(t)\}_{i=1,\dots,N}$ are i.i.d. for $k\in\N$.
\end{lemma}

We point out that this is not the case for the index $i$, both for the particles $u_{ik}^{\epsilon}$ and $\bar{u}^{\epsilon}_{ik}=\Bar{u}^\epsilon_k(x_i,\cdot)$.
Indeed, the laws of $u_{ik}^{\epsilon}$ and $u_{jk}^{\epsilon}$, or $\bar{u}_{ik}^{\epsilon}$ and $\bar{u}_{ik}^{\epsilon}$ respectively, might differ as a result of the $x$ dependence of their defining equations.
Furthermore, even if the points $x_i,x_j\in Q$ are far from each other, namely if $|x_i-x_j|>2\epsilon$, so that their noise terms $W^{\epsilon}_k(x_i,t)$ and $W^{\epsilon}_k(x_j,t)$ are independent, these particles might still be correlated as a result of their initial data.
In fact, from the point of view of modelling in neuroscience, we expect $u_k(x,0)$ to be close to $u_k(y,0)$ for $x$ close to $y$.
This is essentially the issue we pointed out in Remark \ref{remark/n mean field problem}.

We have the following uniform moment estimates both for the true particles $u_{ik}^\epsilon$ and for the McKean--Vlasov particles $\bar{u}_{ik}^\epsilon$.
The case $p=2$ is proved in \cite[Theorem 2.3-4]{carrillo_clini_solem_the_mean_field_limit} and the case of a general $p\geq1$ is proved almost identically.

\begin{proposition}[Uniform moment estimates]
\label{proposition/ uniform moment estimates}
In the Setting \ref{setting} above, for every $T\geq0$ and every $p\geq1$, we have
\begin{equation}
\label{formula/ uniform moment estimates}
    \sup_{i,k}\E\left[\sup_{t\in[0,T]}|u^\epsilon_{ik}|^p\right]+
    \sup_{x\in Q}\E\left[\sup_{t\in[0,T]}|\Bar{u}_k^\epsilon(x,t)|^p\right]
    \leq
    C\left(1+\sup_{x\in Q}\E[|u(x,0)|^p]\right),
\end{equation}
for a constant $C = C(T,p,b,\sigma)$ independent of $M,N$ and $\epsilon$.
\end{proposition}

Furthermore, we have the following estimates for the error between the true and the McKean--Vlasov particles. 
The case $p=2$ is proved in \cite[Theorem 2.6]{carrillo_clini_solem_the_mean_field_limit} and the case $p=4$ is a straightforward adaptation.

\begin{theorem}[Mean squared error estimates for actual particles vs. McKean--Vlasov particles]\label{theorem/Mean squared error estimates for actual particles vs McKean--Vlasov particles}
In the Setting \ref{setting} above, for $p=2,4$ and for any $T>0$, we have
\begin{equation}\label{formula/ mean squared error estimates for actual particles vs McKean--Vlasov particles}
    \sup_{i,k}\E\left[\,\sup_{t\in[0,T]}|u_{ik}^{\epsilon}(t)-\Bar{u}_{ik}^{\epsilon}(t)|^{p}\right]^{\color{black}\nicefrac{1}{p}}\leq C \,\left(\frac{1}{\ngrid^{\frac{\alpha}{d}}}+\frac{1}{\ndelta^{\frac{1}{2}}}\right)\left(1+\sup_{x\in Q}\E\left[|u_k(x,0)|^p\right]^{\nicefrac{1}{p}}\right),
\end{equation}
for a constant $C = C(T,b,\sigma,[u(\cdot,0)]_{\alpha})$, where $[u(\cdot,0)]_{\alpha}$ denotes the H\"older seminorm of $u(\cdot,0)\in C^{\alpha}\big(Q;L^{p}(\Omega)\big)$.
\end{theorem}

We now consider empirical and actual measures associated to the true and McKean--Vlasov particles.
In \cite{carrillo_clini_solem_the_mean_field_limit} it is established that the solution $f(x,t,du)=\law_{\R^{\dirnum}}(\Bar{u}^{\epsilon}(x,t))$ of the Fokker--Planck PDE \eqref{abstract model nonlinear fokker--planck} is independent of $\epsilon$ and it is $\alpha$-H\"older continuous as a function of $x$ with values in $\pr_2\left(C([0,T];\R^{\dirnum})\right)$ endowed with the $2$-Wasserstein metric whenever $u(\cdot,0)\in C^{\alpha}(Q;L^2(\Omega))$.
Furthermore it is proved that $f_{MN}\to f$ in Wasserstein distance of order 1, in expectation.

In the following we will need analogous results for the joint law of two McKean--Vlasov particles 
\begin{equation}
    \label{formula/joint law of mckean--vlasov particles}
    f^{2,\epsilon}(x,y,t,du,dv):=\law_{\R^{2\dirnum}}\left(\Bar{u}^\epsilon(x,t),\Bar{u}^\epsilon(y,t)\right),
\end{equation}
which is a measure on $\R^{2\dirnum}$ depending on $x,y\in Q$, $t\in\R^+$ and $\epsilon>0$, inducing a measure on $Q^2\times\R^{2\dirnum}$ via integration as in \eqref{formula/ measure induced via integration}, and for the joint empirical measure associated to the actual particles
\begin{equation}
    \label{formula/joint empirical measure of true particles}
    f^2_{MN}(t,dx,dy,du,dv):=\frac{1}{M N^2}\sum_{i,j=1}^N\sum_{k=1}^M\delta_{(x_i,x_j,u_{ik}^\epsilon,u_{jk}^\epsilon)}\in\pr\left(Q^2\times\R^{2\dirnum}\right).
\end{equation}

\begin{remark}\label{remark/dependence on epsilon}
As anticipated in the introduction, the joint law $f^{2,\epsilon}(x,y,t,du,dv)$ does depend on the correlation radius $\epsilon>0$, in contrast with the case of the single law $f(x,t,du)$.
Indeed, the independence on $\epsilon$ of $f$ is an artefact of the specific type of interaction among particles, which takes the form, highlighting the drift term only,
\begin{align}\label{formula/specific interaction among particles}
\begin{split}    
    d\Bar{u}^\epsilon(x,t)&=\dots+\phi\left(\int_{Q\times\R^{\dirnum}}b_1(x,y,t,\Bar{u}^{\epsilon}(x,t),v)f(y,t,dv)\,dy\right)\,dt+\dots
    \\
    &=\dots+\phi\left(\int_{Q}\E\left[b_1(x,y,t,u,\Bar{u}^\epsilon(y,t))\right]_{u=\Bar{u}^{\epsilon}(x,t)}\,dy\right)\,dt+\dots
\end{split}
\end{align}
That is to say, the particle $\Bar{u}^\epsilon(x,t)$ interacts with the particle $\Bar{u}^\epsilon(y,t)$ only through its law $f(y,t,dv)$.
Therefore the correlation radius of the noise sources $W^{\epsilon}(x,t)$ and $W^{\epsilon}(y,t)$ affecting the particles plays no role since each of these terms behaves individually as a Brownian motion, regardless of the value of $\epsilon$.
Accordingly, the effect of $\epsilon$ is recovered as soon as we consider the \emph{joint} law of $\Bar{u}^\epsilon(x,t)$ and $\Bar{u}^\epsilon(y,t)$.

We mention that other commonly proposed models (see e.g. \cite{Amari1977,kuhen_ridler_large_deviations,faugeras_inglis_stochastic_neural_field_equations} and the review \cite[Section 6]{bressloff_spatiotemporal_dynamics_of_continuum_neural_fields}) consider a different type of interaction, where $\Bar{u}^\epsilon(x,t)$ does interact with $\Bar{u}^\epsilon(y,t)$ directly and not only through its law, so that the dependence on $\epsilon$ is retained already at the level of the single law $\law_{\R^{\dirnum}}(\Bar{u}^\epsilon(x,t))$.
\end{remark}

 We can promote $f$ and $f^{2,\epsilon}$ to measures on $C([0,T];\R^{\dirnum})$ and $C([0,T];\R^{2\dirnum})$ respectively simply by considering $f(x,\cdot)=\law_{C([0,T];\R^{\dirnum})}(\Bar{u}_k(x,\cdot))$ and $f^{2,\epsilon}(x,y,\cdot)=\law_{C([0,T];\R^{2\dirnum})}(\Bar{u}_k(x,\cdot),\Bar{u}_k(y,\cdot))$.
 Then, from \cite[Theorem 2.5]{carrillo_clini_solem_the_mean_field_limit}, or straightforward adaptations of its arguments in the cases of $p=4$ and of the joint measure, we have the following.

\begin{proposition}[H\"older continuity of the single and joint empirical measure]\label{proposition/ holderianity of the single and joint empirical measure}
In the Setting \ref{setting} above, for $p=2,4$, we have $f\in C^{\alpha}\left(Q;\pr_p\left(C([0,T];\R^{\dirnum})\right)\right)$ and $f^{2,\epsilon}\in C^{\alpha}\left(Q^2;\pr_p\left(C([0,T];\R^{2\dirnum})\right)\right)$ with the following estimates
\begin{align}
    \begin{split}
        \w_p&\Big(C([0,T];\R^{\dirnum})\Big)(f(x,\cdot),f(x',\cdot))
        \leq C\,\bigg(1+\sup_{x\in Q}\E[|u(x,0)|^p]^{\nicefrac{1}{p}}\bigg)\,|x-x'|^\alpha\quad\forall x,x'\in Q,
        \\
        \w_p&\Big(C([0,T];\R^{2\dirnum})\Big)(f(x,y,\cdot),f(x',y',\cdot))
        \\
        &\leq C\left(1+\sup_{x\in Q}\E[|u(x,0)|^p]^{\nicefrac{1}{p}}\right)\left(|x-x'|^\alpha+|y-y'|^{\alpha}+\frac{|x-x'|+|y-y'|}{\epsilon}\right)\quad\forall x,x',y,y'\in Q ,
    \end{split}
\end{align}
for a constant $C = C(T,\rho,b,\sigma,[u(\cdot,0)]_{\alpha})$, where $[u(\cdot,0)]_{\alpha}$ denotes the H\"older seminorm of $u(\cdot,0)\in C^{\alpha}\big(Q;L^{p}(\Omega)\big)$. 
\end{proposition}

Finally, we establish the convergence for the empirical measure. 
We recall that a space dependent measure $g(x,du)\in\pr(\R^{\dirnum})$ induces a measure on $Q\times\R^{\dirnum}$ via formula \eqref{formula/ measure induced via integration} and similarly for a measure $g(x,y,du,dv)\in\pr(\R^{2\dirnum})$ inducing a measure on $Q^2\times\R^{2\dirnum}$.
A straightforward adaptation of \cite[Theorem 2.7]{carrillo_clini_solem_the_mean_field_limit} proves the following.

\begin{theorem}[Convergence of the single and joint empirical measure]\label{theorem/ convergence of the single and joint empirical measure}
In the Setting \ref{setting} above, we have $f_{MN}\to f$ and $f^2_{MN}\to f^{2,\epsilon}$ in Wasserstein distance $\w_2$ in expectation with convergence rates, for every $T\geq0$ and every $q>2$, for a constant $C=C(T,q,\rho,b,\sigma,[u(\cdot,0)]_{\alpha})$,
{\small
\begin{align}\label{formula/ convergence rates in wasserstein distance}
\begin{split}
    \sup_{t\in[0,T]}\E\bigg[\w_2&\left(Q\times\R^{\dirnum}\right)(f_{MN}(t),f(t))^2\bigg]
    \\
    \leq\,
    &
    C\, \left(1+\sup_{x\in Q}\E[|u(x,0)|^2]\right)\, \left(\frac{1}{\sqrt{M}}+\frac{1}{N^{\nicefrac{\alpha}{d}}}\right)^2
    \\
    &
    +C\, \left(1+\sup_{x\in Q}\E[|u(x,0)|^q]^{\nicefrac{2}{q}}\right)\, 
    \begin{cases}
    \displaystyle 
    M^{-\nicefrac{1}{2}}+M^{\nicefrac{-(q-2)}{q}} \quad\text{if $2>\nicefrac{\dirnum}{2}$ and $q\neq 4$,}
    \\
    M^{-\nicefrac{1}{2}}\log(1+M)+M^{\nicefrac{-(q-2)}{q}} \quad\text{if $2=\nicefrac{\dirnum}{2}$ and $q\neq 4$,}
    \\
    M^{-\nicefrac{2}{\dirnum}}+M^{\nicefrac{-(q-2)}{q}} \quad\text{if $2\in(0,\nicefrac{\dirnum}{2})$ and $q\neq \frac{\dirnum}{\dirnum-2}$,}
    \end{cases}
    \\
    \sup_{t\in[0,T]}\E\bigg[\w_2&\left(Q^2\times\R^{2\dirnum}\right)(f_{MN}^2(t),f^{2,\epsilon}(t))^2\bigg]
    \\
    \leq\,
    &
    C\, \left(1+\sup_{x\in Q}\E[|u(x,0)|^2]\right)\, \left(\frac{1}{\sqrt{M}}+\frac{1}{N^{\nicefrac{\alpha}{d}}}+\frac{1}{\epsilon N^{\nicefrac{1}{d}}}\right)^2
    \\
    &
    +C\, \left(1+\sup_{x\in Q}\E[|u(x,0)|^q]^{\nicefrac{2}{q}}\right)\, 
    \begin{cases}
    \displaystyle 
    M^{-\nicefrac{1}{2}}+M^{\nicefrac{-(q-2)}{q}} \quad\text{if $2>\dirnum$ and $q\neq 4$,}
    \\
    M^{-\nicefrac{1}{2}}\log(1+M)+M^{\nicefrac{-(q-2)}{q}} \quad\text{if $2=\dirnum$ and $q\neq 4$,}
    \\
    M^{-\nicefrac{1}{\dirnum}}+M^{\nicefrac{-(q-2)}{q}} \quad\text{if $2\in(0,\dirnum)$ and $q\neq \frac{\dirnum}{\dirnum-1}$.}
    \end{cases}
\end{split}
\end{align}
}
\end{theorem}



\section{The central limit theorem}
\label{section/the central limit theorem}

In this section we prove our main result following the strategy sketched in the introduction.
We will \emph{always} consider the Setting \ref{setting}.

We want to find an equation satisfied by the rescaled fluctuations $\eta_t^{MN}=\cmn(f_{MN}-f)$, for some scaling parameter $\cmn$, and pass this to the limit as $M,N\to\infty$.
We anticipate that the right scaling for the fluctuations is $\cmn=\sqrt{M}$ along a suitable scaling regime $M,N\to\infty$, namely $\sqrt{M}N^{-\nicefrac{\alpha}{d}}\to0$.
However, we choose to keep the scaling factor $\cmn$ implicit so as to highlight why this is the only possible choice in order to obatin a nontrivial limit for the fluctuations (cf. Remark \ref{remark/ rescaling factor sqrt M} below).

Let $\psi:Q\times\R^{\dirn}\to\R$ be sufficiently smooth with $\nabla_u\psi(x,u)\cdot \textbf{n}_{\partial\R^{\dirn}_+}\equiv 0$, for example $\psi\in H_x^{\e}\h^{\kone,\tone}$.
Testing $f_{MN}$ against $\psi$, applying It\^o formula and using the boundary conditions on $\psi$ we obtain
\begin{align}
\label{formula/ weak formulation of empirical measure}
    \begin{split}
        \langle f_{MN}(t),\psi\rangle
        =&
        \frac{1}{MN}\sum_{i=1}^N\sum_{k=1}^M\psi(x_i,u_{ik}^\epsilon(0))
        \\
        &+
        \int_0^t\frac{1}{MN}\sum_{i=1}^N\sum_{k=1}^M\nabla_u\psi(x_i,u_{ik}^\epsilon(r))\,b(x_i,u_{ik}^\epsilon,r,f_{MN})\,dr
        \\
        &+
        \int_0^t\frac{1}{MN}\sum_{i=1}^N\sum_{k=1}^M\Delta_u\psi(x_i,u_{ik}^\epsilon(r))\,\sigma(x_i,u_{ik}^\epsilon,r,f_{MN})^2\,\,dr
        \\
        &+
        \frac{1}{MN}\sum_{i=1}^N\sum_{k=1}^M\int_0^t\nabla_u\psi(x_i,u_{ik}^\epsilon(r))\,\sigma(x_i,u_{ik}^\epsilon,r,f_{MN})\,\,dW^{\epsilon}_{ik}(r)
        \\
        &+
        \frac{1}{MN}\sum_{i=1}^N\sum_{k=1}^M\int_0^t\nabla_u\psi(x_i,u_{ik}^\epsilon(r))\cdot\textbf{n}_{\partial\R_+^{\dirn}}\,\,d\ell^{\epsilon}_{ik}(r)
        \\
        =&
        \langle f_{MN}(0),\psi\rangle + \int_0^t\langle f_{MN}(r),L_r(f_{MN}(r))[\psi]\rangle\, dr +\frac{1}{\cmn}M_t^{MN}(\psi),
    \end{split}
\end{align}
for the martingale term
\begin{equation}
    \label{formula/ real martingale term explicit expression / clt section}
    M_t^{MN}(\psi)=\frac{\cmn}{MN}\sum_{i=i}^N\sum_{k=1}^M\int_0^t\nabla_u\psi(x_i,u_{ik}(r))\,\sigma(x_i,r, u_{ik}(r),f_{MN})\, dW_{ik}(r),
\end{equation}
and the differential operator $L_r(f_{MN}(r))$, depending on $f_{MN}(r,dx,du)\in\pr(Q\times\R^{\dirn})$ and on the time $r$, defined for any $\psi:Q\times\R^{\dirn}\to\R$ sufficiently smooth by
\begin{equation}
    \label{formula/ generator of the particle sde}
    L_r(\mu)[\psi](x,u):=\nabla_u\psi(x,u)\,b(x,r,u,\mu)+\Delta_u\psi(x,u)\,\sigma(x,r,u,\mu)^2\quad\forall\mu\in\pr(Q\times\R^{\dirn}).
\end{equation}

Similarly we test $\psi$, which is admissible because of the boundary conditions, against the weak solution $f(t,x,du)$ of the Fokker--Planck PDE \eqref{abstract model nonlinear fokker--planck}, which is regarded as an element of $\pr(Q\times\R^{\dirn})$ via formula \eqref{formula/ measure induced via integration}.
We obtain
\begin{equation}
    \label{formula/ weak formulation of fokker--planck pde}
    \langle f(t),\psi\rangle=\langle f(0), \psi\rangle + \int_0^t\langle f(r,dx,du),L_r(f(r))[\psi]\rangle \, dr,
\end{equation}
where $L_r(f(r))$ is the differential operator \eqref{formula/ generator of the particle sde} evaluated on $f(r,dx,du)\in\pr(Q\times\R^{\dirn})$.

Combining \eqref{formula/ weak formulation of empirical measure} and \eqref{formula/ weak formulation of fokker--planck pde} yields
\begin{equation}
    \label{formula/ weak formulation of fluctuatios}
    \langle \eta_t^{MN},\psi\rangle=\langle \eta_0^{MN}, \psi\rangle + \int_0^t\cmn \left(\langle f_{MN}(r),L_r(f_{MN})[\psi]\rangle-\langle f(r),L_r(f)[\psi]\rangle \right)\, dr+ M_t^{MN}(\psi).
\end{equation}
Finally we rewrite
\begin{equation}
\label{formula/ rewriting of (f,Lf)-(fMN,LfMN)}
    \cmn \left(\langle f_{MN}(r),L_r(f_{MN})[\psi]\rangle-\langle f(r),L_r(f)[\psi]\rangle \right)
    =
    \langle\eta_r^{MN}(dy,dv),\lin_r(f_{MN}(r),f(r))[\psi](y,v)\rangle
\end{equation}
for the linear differential operator $\lin_t(\mu,\nu)$, depending on $\mu,\nu\in\pr(Q\times\R^{\dirn})$, defined by
\begin{align}
\label{formula/ linearized operator}
\small
    \begin{split}
        \lin_t&(\mu,\nu)[\psi](y,v)
        \\ 
        =&\,\, L_t(\mu)[\psi](y,v)
        \\
        &+ \left\langle\nu(dx,du),\nabla_u\psi(x,u)b_1(x,y,t,u,v)\bigg(\int_0^1\Dot{\phi}\left((1-\lambda)b_1(x,t,u,\nu)+\lambda b_1(x,t,u,\mu)\right)d\lambda\bigg)\right\rangle
        \\
        &+
        \left\langle\nu(dx,du),\Delta_u\psi(x,u)\sigma_1(x,y,t,u,v)2\sigma_0(x,t,u)\bigg(\int_0^1\Dot{\phi}\left((1-\lambda)\sigma_1(x,t,u,\nu)+\lambda \sigma_1(x,t,u,\mu)\right)d\lambda\bigg)\right\rangle
        \\
        &+
        \bigg\langle\!\nu(dx,du),\Delta_u\psi(x,u)\sigma_1(x,y,t,u,v)
        \\
        &\qquad\qquad\quad\quad\cdot\bigg(\!\!\int_0^1\!\!\Dot{\phi}\left((1-\lambda)\sigma_1(x,t,u,\nu)+\lambda \sigma_1(x,t,u,\mu)\right)d\lambda\bigg)\big(\phi(\sigma_1(x,t,u,\nu))+\phi(\sigma_1(x,t,u,\mu))\big)\!\bigg\rangle,
    \end{split}
\end{align}
where we used the shorthand \eqref{formula/ notation for b1(x,t,u,mu)} for $b_1(x,t,u,\mu)$ and the analogous terms.
Equation \eqref{formula/ rewriting of (f,Lf)-(fMN,LfMN)} is justified by direct computation simply by plugging in the expression \eqref{formula/ linearized operator} for $\lin(f_{MN},f)$, which is essentially the linearization of $L(f_{MN})-L(f)$ around $f$.
Combining \eqref{formula/ weak formulation of fluctuatios} and \eqref{formula/ rewriting of (f,Lf)-(fMN,LfMN)} we have proved the following.

\begin{proposition}
\label{porposition/ real semimartingale rewriting of etaMN}
For any $\psi\in H_x^{\e}\h^{\kone,\tone}$ the process $\langle\eta_t^{MN},\psi\rangle$ is a real valued continuous semimartingale with decomposition
\begin{equation}
    \label{formula/ real semimartingale rewriting of etaMN}
    \langle \eta_t^{MN},\psi\rangle=\langle \eta_0^{MN}, \psi\rangle + \int_0^t\langle\eta_r^{MN},\lin_r(f_{MN},f)[\psi]\,\rangle\,\, dr+ M_t^{MN}(\psi).
\end{equation}
\end{proposition}

The next step is to understand the behavior of the martingale term $M_t^{MN}(\psi)$ in the $M,N\to\infty$ limit.
We start by analyzing its quadratic variation.
For any $\varphi,\psi\in H_x^{\e}\h^{\kone,\tone}$, from the expressions \eqref{formula/ real martingale term explicit expression / clt section} and \eqref{formula/ space colored time white noise}, we compute:
\begin{align}
\label{formula/ quadratic variation of real martingale term}
\small
    \begin{split}
        \langle&  M_t^{MN}(\varphi),M_t^{MN}(\psi)\rangle
        \\
        =&
        \frac{\cmn^2}{M}\int_0^t\frac{1}{MN^2}\sum_{i,j=1}^N\sum_{k=1}^M\left(\nabla_u\varphi(x_i,u_{ik}^\epsilon)\sigma(x_i,r,u_{ik}^\epsilon,f_{MN})\right)
        \\
        &\qquad\quad\qquad\qquad\qquad\qquad\cdot\left(\nabla_u\psi(x_j,u_{jk}^\epsilon)\sigma(x_j,r,u_{jk}^\epsilon,f_{MN})\right)\Big(\int_{\R^d}\rho(z+\epsilon^{-1}(x_i-x_j))\rho(z)\,dz \, C_{\rho}\Big)\,dr
        \\
        =&
        \frac{\cmn^2}{M}\int_0^t\int_{Q^2\times\R^{2\dirn}}\!\!\!\!\!\!\!\!\!\!\!\!\!\!\!\!\!\left(\nabla_u\varphi(x,u)\sigma(x,r,u,f_{MN})\right)\left(\nabla_u\psi(y,v)\sigma(y,r,v,f_{MN})\right)R^{\epsilon}(x-y)\,df_{MN}^2(r,dx,dy,du,dv)\,dr
        \\
        =&
        \frac{\cmn^2}{M}\int_0^t\!\!\!\left\langle\left(\nabla_u\varphi(x,u)\sigma(x,r,u,f_{MN})\right)\left(\nabla_u\psi(y,v)\sigma(y,r,v,f_{MN})\right)R^{\epsilon}(x-y),\,df_{MN}^2(r,dx,dy,du,dv)\right\rangle\,dr,
    \end{split}
\end{align}
where we have defined
\begin{equation}
    \label{formula/ correlation function R epsilon}
    R^{\epsilon}(x):=C_{\rho}\int_{\R^d}\rho(z+\epsilon^{-1}x)\rho(z)\,dz.
\end{equation}
The expression \eqref{formula/ quadratic variation of real martingale term} and the convergence $f_{MN}\to f$ and $f_{MN}^2\to f^{2,\epsilon}$ from Theorem \ref{theorem/ convergence of the single and joint empirical measure} suggest that the quadratic variation $\langle M_t^{MN}(\varphi),M_t^{MN}(\psi)\rangle$ converges to the \emph{deterministic} finite variation function
\begin{align}\label{formula/ quadratic variation of real gaussian limit}
\small
\begin{split}
    g_t^\epsilon(\varphi,\psi):=
    &\int_0^t\int_{Q^2\times\R^{2\dirn}}\!\!\!\!\!\!\!\!\left(\nabla_u\varphi(x,u)\sigma(x,r,u,f)\right)\left(\nabla_u\psi(y,v)\sigma(y,r,v,f)\right)R^{\epsilon}(x-y)\,df^{2,\epsilon}(r,dx,dy,du,dv)\,dr
        \\
    =&
    \int_0^t\left\langle\left(\nabla_u\varphi(x,u)\sigma(x,r,u,f)\right)\left(\nabla_u\psi(y,v)\sigma(y,r,v,f)\right)R^{\epsilon}(x-y),\,df^{2,\epsilon}(r,dx,dy,du,dv)\right\rangle\,dr.
\end{split}
\end{align}
The next proposition establishes this rigorously.

\begin{proposition} \label{proposition/ pointwise convergence of the quad variation}
For any $\varphi,\psi\in H_x^{\e}\h^{\ktwo,\ttwo}$, for any $T\geq0$ we have
\begin{align}\label{formula/ stronger pointwise convergence of the quad variation}
\small
\begin{split}
    \lim_{M,N\to\infty}\sup_{r\in[0,T]}\!\!\E\bigg[\Big|
    \Big\langle\!\!&\left(\nabla_u\varphi(x,u)\sigma(x,r,u,f_{MN})\right)\left(\nabla_u\psi(y,v)\sigma(y,r,v,f_{MN})\right)R^{\epsilon}(x-y),\,f_{MN}^2(r,dx,dy,du,dv)\Big\rangle
    \\
    &-\Big\langle\!\!\left(\nabla_u\varphi(x,u)\sigma(x,r,u,f)\right)\left(\nabla_u\psi(y,v)\sigma(y,r,v,f)\right)R^{\epsilon}(x-y),\,f^{2,\epsilon}\!(r,dx,dy,du,dv)\Big\rangle
    \Big|\bigg]
    \!=\!
    0.
\end{split}
\end{align}
In particular we have 
\begin{equation}\label{formula/ pointwise convergence of the quad variation}
    \lim_{M,N\to\infty}\E\left[\sup_{t\in[0,T]}\left|\langle M_t^{MN}(\varphi),M_t^{MN}(\psi)\rangle-g^{\epsilon}(\varphi,\psi)(t)\right|\right]=0.
\end{equation}
\end{proposition}

Given an admissible test function $\psi$, standard probability theory ensures there exists a (unique in law) Gaussian martingale $G_t^\epsilon(\psi)$ having the increasing function $g^\epsilon_t(\psi,\psi)$ as its quadratic variation.
Proposition \ref{proposition/ pointwise convergence of the quad variation} hints that $M_t^{MN}(\psi)$ converges to the process $G_t^{\epsilon}(\psi)$.
In turn, assuming that $\eta_t^{MN}$ has some limit $\eta_{t}^\epsilon$, the convergence $f_{MN}\to f$ suggests that we pass to the limit in equation \eqref{formula/ real semimartingale rewriting of etaMN} and obtain the following expression
\begin{equation}
    \label{formula/ real semimartingale rewriting of eta epsilon}
    \langle\eta_t^\epsilon,\psi\rangle = \langle\eta_0,\psi\rangle + \int_0^t\langle \eta_r^\epsilon, \lin_r(f(r),f(r))[\psi]\rangle \,dr + G_t^\epsilon(\psi).
\end{equation}
The remaining part of this section is devoted to making this argument rigorous.
The missing ingredient is essentially the tightness of the processes involved.

First of all, we pass everything to the level of distributions.
We start with the martingale term.
Indeed, it is not hard to see that, for a fixed time $t$, the association $H_x^{\e}\h^{\kone,\tone}\ni\psi\mapsto M_t^{MN}(\psi)$ defines an $H_x^{-\e}\h^{-\kone,\tone}$-valued random variable $M_t^{MN}$.
We have the following result.

\begin{proposition}\label{proposition/ H valued martingale MtMN}
The process $M_t^{MN}$ is an $H_x^{-\e}\h^{-\kone,\tone}$-valued continuous square integrable martingale.
Its quadratic variation is the $H_x^{-\e}\h^{-\kone,\tone}\otimes H_x^{-\e}\h^{-\kone,\tone}$-valued process $\langle M_t^{MN}\rangle$ defined by
\begin{equation}
    \label{formula/ quadratic variation of H valued martingale term}
    \langle M_t^{MN}\rangle(\varphi,\psi):=\langle M_t^{MN}(\varphi),M_t^{MN}(\psi)\rangle.
\end{equation}
For every $T\geq0$, for a constant $C=C(T,b,\sigma)$, we have
\begin{align}\label{formula/ uniform estimate on martingale term}
\begin{split}
\sup_{M,N}\E\!\left[\!\sup_{t\in[0,T]}\!\!\|M_t^{MN}\|_{H_x^{-\e}\h^{-\kone,\tone}}^2\!\right]\!
\leq\, C
\sup_{M,N}\frac{\cmn^2}{M}\bigg(1\!+\!\sup_{x\in Q}\E\big[|u(x,0)|^{2\tone+2}\big]\bigg).
\end{split}
\end{align}
\end{proposition}

Next we want to study the fluctuation process $\eta_t^{MN}\in\pr(Q\times\R^{\dirn})\subseteq H_x^{-\e}\h^{-\kone,\tone}$, obtain bounds on its norm and pass \eqref{formula/ real semimartingale rewriting of etaMN} to an equation at the distributional level.
This will be done in two steps, working between the nested spaces $H_x^{\eplus}\h^{\ktwo,\ttwo}\subseteq H_x^{\e}\h^{\kone,\tone}$.

\begin{lemma}
    \label{proposition/ first bounds on etaMN}
For every $T\geq0$, for a constant $C=C(T,b,\sigma)$, we have
\begin{equation}
    \label{formula/ first bounds on etaMN}
    \sup_{MN}\sup_{t\in[0,T]}\E\left[\|\eta_t^{MN}\|^2_{H^{-\e}\h^{-\kone,\tone}}\right]
    \leq \,C\,
    \sup_{M,N}\left(\frac{\cmn^2}{M}+\frac{\cmn^2}{N^{\nicefrac{2\alpha}{d}}}\right)\bigg(1+\sup_{x\in Q}\E\big[|u(x,0)|^{4\tone}\big]\bigg).
\end{equation}
\end{lemma}

\begin{remark}\label{remark/ rescaling factor sqrt M}
Lemma \ref{proposition/ first bounds on etaMN} implies that, in order to keep the fluctuations bounded in a suitable distribution space, we have to impose the constraint $\displaystyle \limsup_{M,N\to\infty}\left(\frac{\cmn^2}{M}+\frac{\cmn^2}{N^{\nicefrac{2\alpha}{d}}}\right)<\infty$.
On the other hand, estimate \eqref{formula/ uniform estimate on martingale term} implies that, in order for the martingale term not to vanish, which would imply $\eta_{\cdot}^{MN}\to0$ and result in a suboptimal expansion \eqref{formula/ first order expansion} for the empirical measure, we have to impose $\cmn\gtrsim\sqrt{M}$.
The requirements on the scaling regime are therefore $\cmn\simeq\sqrt{M}$ and $\displaystyle \limsup_{M,N\to\infty}\sqrt{M}N^{-\nicefrac{\alpha}{d}}<\infty$.
In fact, we will see that the convergence of the initial data $\eta_0^{MN}$ further requires $\displaystyle \lim_{M,N\to\infty}\sqrt{M}N^{-\nicefrac{\alpha}{d}}=0$.

We can make sense of these constraints as follows.
As explained in Remark \ref{remark/n mean field problem}, our situation essentially corresponds to $N$ classical mean field problems with $M$ particles each, one for each column of $M$ neurons at location $x_i$, interacting among themselves.
The decay rate $\sqrt{M}$ is the usual decay rate for a standard mean field setting with $M$ particles subjected to independent noise sources $W_k^{\epsilon}(x_i,t)$ for $k=1,\dots,M$.
Therefore, if we want to see a nontrivial behavior of the rescaled fluctuations $\cmn(f_{MN}-f)$ and also to keep these bounded, we have to choose $\cmn=\sqrt{M}$.

However, we have $N$ such mean field clusters interacting with each other and two different columns at nearby locations $x_i,x_j$ are initiated with possibly correlated initial data and sense correlated sources of noise $\{W_k^{\epsilon}(x_i,t)\}_{k=1\dots M}$ and $\{W_k^{\epsilon}(x_j,t)\}_{k=1\dots M}$ with correlation strength proportional to $|x_i-x_j|^{\alpha}\simeq N^{-\nicefrac{\alpha}{d}}$.
The condition $\limsup \sqrt{M}N^{-\nicefrac{\alpha}{d}}<\infty$ serves therefore to prevent these correlated random clusters from interacting constructively and making the $\sqrt{M}$-rescaled fluctuations yet unbounded.
\end{remark}

We can obtain better estimates for the process $\eta_t^{MN}$ by exploiting the semimartingale decomposition \eqref{formula/ real semimartingale rewriting of etaMN}.
For this, we first need the following lemma on the norm of the operators $\lin_t(\mu,\nu)$ defined in \eqref{formula/ linearized operator}.

\begin{lemma}\label{proposition/ norm bounds on the linearized operator}
The linear operators $\lin_t(f_{MN},f)$ and $\lin_t(f,f)$ define two maps $H_x^{e}\h^{\ktwo,\ttwo}\to H_x^{\e}\h^{\kone,\tone}$ and two maps $H_x^{e}\h^{\ktwo+2,\ttwo-2}\to H_x^{\e}\h^{\ktwo,\ttwo}$ satisfying the bounds,
{\small
\begin{align}\label{formula/ norm of linearized operators 1}
\begin{split}
\sup_{t\in[0,T]}\sup_{M,N}\underset{\omega\in\Omega}{\esssup}\|\lin_t(f_{MN},f)\|_{H_x^{\e}\h^{\ktwo,\ttwo},H_x^{\e}\h^{\kone,\tone}} + &\sup_{t\in[0,T]}\|\lin(f,f)\|_{H_x^{\e}\h^{\ktwo,\ttwo},H_x^{\e}\h^{\kone,\tone}}
\\
&\leq C\bigg(1+\sup_{x\in Q}\E\left[|u(x,0)|^{4+2\tone}\right]\bigg)^{\nicefrac{1}{2}},
\\
\sup_{t\in[0,T]}\sup_{M,N}\underset{\omega\in\Omega}{\esssup}\|\lin_t(f_{MN},f)\|_{H_x^{\e}\h^{\ktwo+2,\ttwo-2},H_x^{\e}\h^{\ktwo,\ttwo}} +& \sup_{t\in[0,T]}\!\!\|\lin(f,f)\|_{H_x^{\e}\h^{\ktwo+2,\ttwo-2},H_x^{\e}\h^{\ktwo,\ttwo}}
\\
&\leq C\bigg(1+\sup_{x\in Q}\E\left[|u(x,0)|^{4+2\tone}\right]\bigg)^{\nicefrac{1}{2}},
\end{split}
\end{align}
}
for a constant ${\small C=C(T,b,\sigma,\dirnum,d,\alpha)}$ independent of the event $\omega\in\Omega$.
\end{lemma}

Using these bounds, we now obtain a semimartingale expression for $\eta_t^{MN}$ in the bigger space $H_x^{-\e}\h^{-\ktwo,\ttwo}\supseteq H_x^{-\e}\h^{-\kone,\tone}$.

\begin{proposition}\label{proposition/ H valued semimartingale etaMN}
The process $\eta_t^{MN}$ is an $H_x^{-\e}\h^{-\ktwo,\ttwo}$-valued continuous square integrable semimartingale with decomposition
\begin{equation}
    \label{formula/ H valued semimartingale rewriting of etaMN}
    \eta_t^{MN}=\eta_0^{MN}+\int_0^t\lin_r(f_{MN},f)^*[\eta^{MN}_r]\,dr + M^{MN}_t,
\end{equation}
where $\lin_r(f_{MN},f)^*:H_x^{-\e}\h^{-\kone,\tone}\to H_x^{-\e}\h^{-\ktwo,\ttwo}$ is the adjoint of the operator $\lin_r(f_{MN},f)$.
For every $T\geq0$, for a constant $C=C(T,b,\sigma,d,\dirnum,\alpha)$, we have
\begin{equation}\label{formula/ second bounds on etaMN}
\small
    \sup_{MN}\E\bigg[\sup_{t\in[0,T]}\|\eta_t^{MN}\|^2_{H^{-\e}\h^{-\ktwo,\ttwo}}\bigg]
    \leq C
    \sup_{M,N}\left(\frac{\cmn^2}{M}+\frac{\cmn^2}{N^{\nicefrac{2\alpha}{d}}}\right)\bigg(1+\sup_{x\in Q}\E\big[|u(x,0)|^{4\tone}\big]\bigg).
\end{equation}
\end{proposition}

In order to pass to the limit in the expression \eqref{formula/ H valued semimartingale rewriting of etaMN} we exploit the tightness of the terms involved.
This is proved with the previous estimates, the compactness of the embedding $H_x^{-\e}\h^{-\kone,\tone}\subseteq H_x^{-\e}\h^{-\ktwo,\ttwo}$ and the Aldous criterion (see e.g. \cite{Joffe_metivier_weak_convergence_of_sequence_of_semimartingales}).

\begin{proposition} \label{proposition/ tightness of the processes}
For $M,N\in\N$ the random variables $\small M_t^{MN}\in C\big([0,\infty];H_x^{-e_{\alpha}^+}\h^{-\ktwo,\ttwo}\big)$ are tight.
Furthermore, the random variables ${\small\eta_t^{MN}\in C\big([0,\infty];H_x^{-e_{\alpha}^+}\h^{-\ktwo,\ttwo}\big)}$ are tight along a scaling regime such that ${\small\displaystyle\limsup_{M,N\to\infty}\sqrt{M}N^{\nicefrac{-\alpha}{d}}<\infty}$.
\end{proposition}

Next we use Prokhorov's theorem and the convergence \eqref{formula/ pointwise convergence of the quad variation} of the quadratic variation $\langle M_t^{MN}\rangle(\varphi,\psi)$ to establish that $M_t^{MN}$ does have a unique limit in $C\big([0,T];H_x^{-\eplus}\h^{-\ktwo,\ttwo}\big)$.
The first step is to identify uniquely this limit.
The following lemma is proved similarly to estimate \eqref{formula/ uniform estimate on martingale term}.

\begin{lemma}\label{proposition/ quadratic variation of the gaussian limit}
The association $H^{\e}_x\h^{\kone,\tone}\ni\varphi,\psi\mapsto g^{\epsilon}_t(\varphi,\psi)$ defines a \emph{deterministic} increasing positive definite function $g_t^\epsilon$ with values in $H^{-\e}_x\h^{-\kone,\tone}\otimes H^{-\e}_x\h^{-\kone,\tone}$.
\end{lemma}

Thus $g_t^\epsilon$ is an admissible covariation function in $H^{-\e}_x\h^{-\kone,\tone}$ and standard probability theory (see e.g. \cite[Chapter 3]{daprato_zabczyk_1992}) ensures there exists a unique-in-law $H^{-\e}_x\h^{-\kone,\tone}$-valued Gaussian martingale $G_t^\epsilon$ with quadratic variation $g_t^\epsilon$.

The uniform bound \eqref{formula/ uniform estimate on martingale term}, Proposition \ref{proposition/ tightness of the processes} and Prokhorov's theorem imply that, along subsequences, the process $M_t^{MN}$ converges both in law in $C\big([0,T];H_x^{-\eplus}\h^{-\ktwo,\ttwo}\big)$ and weakly in the Hilbert space $\mathcal{M}_T^2\big(H_x^{\e}\h^{-\kone,\tone}\big)$ of square integrable continuous martingales to (possibly different) limit martingales.
In addition, Proposition \ref{proposition/ pointwise convergence of the quad variation} implies that along any such subsequence  we have $\langle M_t^{MN}\rangle(\varphi,\psi)\to g_t^\epsilon(\varphi,\psi)$ for every $\varphi,\psi\in H^{\e}_x\h^{\ktwo,\ttwo}$.
Hence, along any subsequence, the limit martingale must be $G_t^\epsilon$ and a sub-subsequence argument readily proves that the whole sequence $M_t^{MN}$ is converging to $G_t^\epsilon$.
Hence we have proved the following.

\begin{proposition}\label{proposition/ clt for the martingale}
We have $\lim_{M,N\to\infty} M_t^{MN}=G_t^\epsilon$ in law in $C\big([0,T];H_x^{-\eplus}\h^{-\ktwo,\ttwo}\big)$, with no constraints on the scaling regime for $M,N\to\infty$.
\end{proposition}

Finally we exploit the tightness of $\eta_t^{MN}$ and the convergence $f_{MN}\to f$ and $M_t^{MN}\to G_t^\epsilon$ to pass to the limit in equation \eqref{formula/ H valued semimartingale rewriting of etaMN} and show that the fluctuations $\eta_t^{MN}$ have a unique limit.
We first establish the convergence of the initial data $\eta_0^{MN}$.
This is obtained with adaptations of the standard argument for Lévy's classical central limit theorem. 

\begin{lemma}\label{proposition/ clt for the initial data}
Consider a scaling regime $M,N\to\infty$ such that $\lim \sqrt{M}N^{-\nicefrac{\alpha}{d}}=0$, then we have
\begin{equation}
    \label{formula/ clt for the initial data}
    \eta_0^{MN}\to\eta_0\quad\text{in law in } H_x^{-\eplus}\h^{-\ktwo,\ttwo},
\end{equation}
where $\eta_0\sim \mathcal{N}(0,\mathcal{Q})$ is a Gaussian r.v. in $H_x^{-\e}\h^{-\kone,\tone}$ with mean zero and covariance 
{\small
\begin{align}\label{formula/ clt covariance of the limit initial data}
\small
\begin{split}
    \mathcal{Q}(\varphi,\psi)&=\E\bigg[\Big(\int_Q\varphi(x,u_k(x,0))-\E[\varphi(x,u_k(x,0))]\,dx\Big)\Big(\int_Q\psi(y,u_k(y,0))-\E[\psi(y,u_k(y,0))]\,dy\Big)\bigg]
    \\
    &=
    \!\int_{Q^2\times\R^{2\dirnum}}\!\!\!\bigg(\!\varphi(x,u)\!-\!\!\int_{\R^{\dirnum}}\!\!\!\!\!\varphi(x,w)f_0(x,dw)\bigg)\bigg(\!\psi(y,v)\!-\!\int_{\R^{\dirnum}}\!\!\!\!\!\psi(y,w)f_0(y,dw)\bigg) f_0^2(x,y,du,dv)\,dx\,dy.
\end{split}
\end{align}
}
In particular $\eta_0$ is independent of $\epsilon$ since both $f_0(x,du)=\law_{\R^{\dirnum}}(u_k(x,0))$ and $f_0^2(x,y,du,dv)=\law_{\R^{2\dirnum}}(u_k(x,0),u_k(y,0))$ are.
\end{lemma}

Finally we establish our main result: a central limit theorem for the whole fluctuation process $\eta_t^{MN}$.
Here below, the well-posedness of equation \eqref{formula/ langevin spde / clt section} is guaranteed by classical stochastic analysis in infinite dimension (see e.g. \cite[Chapter 5]{daprato_zabczyk_1992}).

\begin{theorem} \label{theorem/ clt/ clt section}
Consider the Setting \ref{setting}.
For every $\epsilon>0$ fixed, as ${\small M,N\to\infty}$ with scaling regime ${\small \lim \sqrt{M}N^{-\nicefrac{\alpha}{d}}=0}$, we have
\begin{equation}
    \label{formula/ clt for fluctuation process}
    \eta_t^{MN}\to\eta_t^\epsilon\quad\text{in law in } C\big([0,T];H_x^{-\eplus}\h^{-\ktwo,\ttwo}\big),
\end{equation}
where $\eta^\epsilon$ is the unique weak solution in $H_x^{-\e}\h^{-(\ktwo+2),\ttwo-2}$ of the Langevin SPDE 
\begin{align}\label{formula/ langevin spde / clt section}
    \eta_t^{\epsilon}=\eta_0+\int_0^t\lin_r(f,f)^*[\eta_r^{\epsilon}]\,dr + G_t^\epsilon.
\end{align} 
\end{theorem}



\section{Proof of the results}
\label{section/proof of the results}

\begin{proof}[\textbf{Proof of Proposition \ref{proposition/ pointwise convergence of the quad variation}}]
\noindent
By polarization it is sufficient to prove the result when $\varphi=\psi$.
We introduce the following random time dependent probability measures on $Q^2\times\R^{2\dirnum}$:
\begin{align}\label{pointwise convergence of the quad variation/ proof 1}
    \begin{split}
        \bar{f}^2_{MN}(t)=\frac{1}{MN^2}\sum_{i,j=1}^N\sum_{k=1}^M\delta_{\big(x_i,x_j,\Bar{u}^{\epsilon}_{ik}(t),\Bar{u}_{jk}^{\epsilon}(t)\big)}\quad\text{ and } \quad\bar{f}^2_{N}(t)=\frac{1}{N^2}\sum_{i,j=1}^N\delta_{(x_i,x_j)}\otimes f^{2,\epsilon}(x_i,x_j,t,du,dv).
    \end{split}
\end{align}
We consider the splitting 
\begin{align}\label{pointwise convergence of the quad variation/ proof 2}
\small
\begin{split}
    \Big\langle\!\!(\nabla_u&\psi(x,u)\sigma(x,r,u,f_{MN}))\left(\nabla_u\psi(y,v)\sigma(y,r,v,f_{MN})\right)R^{\epsilon}(x-y),\,df_{MN}^2(r,dx,dy,du,dv)\Big\rangle
    \\
    &-\Big\langle\!\!\left(\nabla_u\psi(x,u)\sigma(x,r,u,f)\right)\left(\nabla_u\psi(y,v)\sigma(y,r,v,f)\right)R^{\epsilon}(x-y),\,f^{2,\epsilon}\!(r,dx,dy,du,dv)\Big\rangle
    \\
    =&
    \Big\langle\!\!\left(\nabla_u\psi(x,u)\sigma(x,r,u,f_{MN})\right)\left(\nabla_u\psi(y,v)(\sigma(y,r,v,f_{MN})-\sigma(y,r,v,f))\right)R^{\epsilon}(x-y),\,f_{MN}^2\Big\rangle
    \\
    &+
    \Big\langle\!\!\left(\nabla_u\psi(x,u)(\sigma(x,r,u,f_{MN})-\sigma(x,r,u,f))\right)\left(\nabla_u\psi(y,v)\sigma(y,r,v,f)\right)R^{\epsilon}(x-y),\,f_{MN}^2\Big\rangle
    \\
    &+
    \Big\langle\!\!\left(\nabla_u\psi(x,u)\sigma(x,r,u,f)\right)\left(\nabla_u\psi(y,v)\sigma(y,r,v,f)\right)R^{\epsilon}(x-y),\,f_{MN}^2-\Bar{f}_{MN}^2\Big\rangle
    \\
    &+
    \Big\langle\!\!\left(\nabla_u\psi(x,u)\sigma(x,r,u,f)\right)\left(\nabla_u\psi(y,v)\sigma(y,r,v,f)\right)R^{\epsilon}(x-y),\,\Bar{f}_{MN}^2-\Bar{f}_N^2\Big\rangle
    \\
    &+
    \Big\langle\!\!\left(\nabla_u\psi(x,u)\sigma(x,r,u,f)\right)\left(\nabla_u\psi(y,v)\sigma(y,r,v,f)\right)R^{\epsilon}(x-y),\,\Bar{f}_N^2-f^{2,\epsilon}\Big\rangle
    \\
    =& E_r^1+E_r^2+E_r^3+E_r^4+E_r^5.
\end{split}
\end{align}
We consider each term $E_r^i$ separately and we show that they satisfy $\sup_{r\in[0,T]}\E\big[|E_r^i|\big]\to0$.

For the first term we compute, for a constant $C=C(T,b,\sigma)$,
\begin{align}\label{pointwise convergence of the quad variation/ proof 3}
\small
\begin{split}
    &\sup_{r\in[0,T]}\E\big[|E_r^1|\big]
    \\
    &\leq
    \sup_{r\in[0,T]}\!\!\E\bigg[ \int_{Q^2\times\R^{2\dirnum}}\!\!\!\!\!\!\!\!\!\!\!\!\!\!\!\!\!\!|\nabla_u\psi(x,u)|\,|\sigma(x,r,u,f_{MN})|\,|\nabla_u\psi(y,v)|\,|\sigma(y,r,v,f_{MN})-\sigma(y,r,v,f))|\,|R^{\epsilon}(x-y)|\,\,df_{MN}^2\bigg]
    \\
    &\leq C\,\sup_{r\in[0,T]}
    \E\bigg[\int_{Q^2\times\R^{2\dirnum}}\|\psi\|_{H_x^{\e}\h^{\ktwo,\ttwo}}^2(1+|u|^{\ttwo+1})(1+|v|^{\ttwo})\big(\w_1(Q\times\R^{\dirnum})(f_{MN},f)\big)^\alpha\,\, df_{MN}^2\bigg]
    \\
    &\leq C\,\|\psi\|_{H_x^{\e}\h^{\ktwo,\ttwo}}^2
    \!\!\!\sup_{r\in[0,T]}\!\E\bigg[\int_{Q^2\times\R^{2\dirnum}}\!\!\!\!\!\!\!\!1+|u|^{4\ttwo+2}+|v|^{4\ttwo}\,\, df_{MN}^2\bigg]^{\nicefrac{1}{2}}
    \!\!\!\sup_{r\in[0,T]}\E\Big[\big(\w_1(Q\times\R^{\dirnum})(f_{MN},f)\big)^2\Big]^{\nicefrac{\alpha}{2}}
    \\
    &\leq C\,\|\psi\|_{H_x^{\e}\h^{\ktwo,\ttwo}}^2
    \bigg(1+\sup_{x\in Q}\E\big[|u_k(x,0)|^{4\ttwo+2}\big]\bigg)^{\nicefrac{1}{2}}
    \sup_{r\in[0,T]}\E\Big[\big(\w_2(Q\times\R^{\dirnum})(f_{MN},f)\big)^2\Big]^{\nicefrac{\alpha}{2}}.
\end{split}
\end{align}
In the second passage we used the embedding \eqref{formula/ continuous embedding of working spaces into bounded functions} to bound $|\nabla_u\psi|$, the linear growth and Lipschitz properties \eqref{formula/ properties of b and sigma} of $\sigma$ and the boundedness of $R^{\epsilon}$.
The third passage follows from several application of H\"older's and Young's inequality.
In the last passage we used the definition of $f_{MN}^2$ and the moment estimates \eqref{formula/ uniform moment estimates}.

Swapping the roles of $(x,u)$ and $(y,v)$ and of $f_{MN}$ and $f$ in \eqref{pointwise convergence of the quad variation/ proof 3}, we obtain an identical estimate for the second term $E_r^2$. 
In particular, we notice that the right hand side of \eqref{pointwise convergence of the quad variation/ proof 3} vanishes as $M,N\to\infty$ thanks to Theorem \ref{theorem/ convergence of the single and joint empirical measure}.

For the third term in \eqref{pointwise convergence of the quad variation/ proof 2} we compute, for a constant $C=C(T,b,\sigma)$,
\begin{align}\label{pointwise convergence of the quad variation/ proof 4}
\small
\begin{split}
    |E_r^3|
    \leq\!
    \frac{C}{MN^2}\!\!\sum_{i,j=1}^N\!\sum_{k=1}^M\!
    \bigg(
    &|\nabla_u\psi(x_i,u_{ik}^{\epsilon})-\nabla_u\psi(x_i,\bar{u}_{ik}^{\epsilon})||\sigma(x_i,u_{ik}^{\epsilon},r,f)\nabla_u\psi(x_j,u_{jk}^{\epsilon}\sigma(x_j,u_{jk}^{\epsilon},r,f)|
    \\
    &+
    |\nabla_u\psi(x_i,\bar{u}_{ik}^{\epsilon})||\sigma(x_i,u_{ik}^{\epsilon},r,f)-\sigma(x_i,\bar{u}_{ik}^{\epsilon},r,f)||\nabla_u\psi(x_j,u_{jk}^{\epsilon})\sigma(x_j,u_{jk}^{\epsilon},r,f)|
    \\
    &+
    |\nabla_u\psi(x_i,\bar{u}_{ik}^{\epsilon})||\sigma(x_i,\bar{u}_{ik}^{\epsilon},r,f)||\nabla_u\psi(x_j,u_{jk}^{\epsilon})-\nabla_u\psi(x_j,\bar{u}_{jk}^{\epsilon})||\sigma(x_j,u_{jk}^{\epsilon},r,f)|
    \\
    &+
    |\nabla_u\psi(x_i,\bar{u}_{ik}^{\epsilon})||\sigma(x_i,\bar{u}_{ik}^{\epsilon},r,f)||\nabla_u\psi(x_j,\bar{u}_{jk}^{\epsilon})||\sigma(x_j,u_{jk}^{\epsilon},r,f)-\sigma(x_j,\Bar{u}_{jk}^{\epsilon},r,f)|
    \bigg)
    \\
    \leq\!
    \frac{C}{MN^2}\!\!\sum_{i,j=1}^N\!\sum_{k=1}^M\!
    \big(1&+\|\psi\|^3_{H_x^{\e}\h^{\ktwo,\ttwo}}\big)\big(1+|u_{ik}^{\epsilon}|+|\Bar{u}_{ik}^{\epsilon}|+|u_{jk}^{\epsilon}|+|\Bar{u}_{jk}^{\epsilon}|\big)^{3\ttwo+1}\big(|u_{ik}^{\epsilon}\!-\!\Bar{u}_{ik}^{\epsilon}|+|u_{jk}^{\epsilon}\!-\!\Bar{u}_{jk}^{\epsilon}|\big).
\end{split}
\end{align}
In the first passage we unfolded the term $E_r^3$ using the definition of $f^2_{MN}$ and $\Bar{f}_{MN}^2$, and then we added and subtracted several mixed term.
In the second passage we used the Lipschitz and linear growth properties \eqref{formula/ properties of b and sigma} of $\sigma$, the definition of the operators $V_{(x,u)}^1$ and $V_{(x,u),(y,v)}^{1,\text{dif}}$ and the estimates \eqref{formula/ norm of evaluation operators} on their norms, and several applications of H\"older's and Young's inequality.

We now take the expectation and the supremum in time of \eqref{pointwise convergence of the quad variation/ proof 4}, and then we apply the Cauchy--Schwarz inequality and the moment estimates \eqref{formula/ uniform moment estimates} to obtain, for a constant $C=C(T,b,\sigma)$,
{\small
\begin{align}\label{pointwise convergence of the quad variation/ proof 5}
\begin{split}
\sup_{r\in[0,T]}\E\big[|E_r^3|\big]\leq C\,\big(1&+\|\psi\|^3_{H_x^{\e}\h^{\ktwo,\ttwo}}\big)\bigg(1+\sup_{x\in Q}\E\big[|u_k(x,0)|^{6\ttwo+2}\big]^{\nicefrac{1}{2}}\bigg)\sup_{i,k}\sup_{r\in[0,T]}\E\big[|u_{ik}^{\epsilon}(r)-\Bar{u}_{ik}^{\epsilon}(r)|^2\big]^{\nicefrac{1}{2}}.
\end{split}
\end{align}
}
We note the right-hand side of the inequality vanishes in the $M,N\to\infty$ limit thanks to Theorem \ref{theorem/Mean squared error estimates for actual particles vs McKean--Vlasov particles}.

For the fourth term in \eqref{pointwise convergence of the quad variation/ proof 2} we compute, for a constant $C=C(T,b,\sigma)$,
{\small
\begin{align}\label{pointwise convergence of the quad variation/ proof 6}
\begin{split}
&\sup_{r\in[0,T]}\E\big[|E_r^4|^2\big]
\\
&\leq C
\frac{1}{N^2}\sum_{i,j=1}^N\sup_{r\in[0,T]}\!\!\!\E\Bigg[
\bigg|\frac{1}{M}\sum_{k=1}^M\bigg(\left(\nabla_u\psi(x,u)\sigma(x,r,u,f)\right)\left(\nabla_u\psi(y,v)\sigma(y,r,v,f)\right)
\\
&\qquad\qquad\quad\qquad\quad\qquad\qquad
-\int_{\R^{2\dirnum}}\!\!\!\!\!\!\!\!\!\left(\nabla_u\psi(x,u)\sigma(x,r,u,f)\right)\left(\nabla_u\psi(y,v)\sigma(y,r,v,f)\right)\,f^{2,\epsilon}(r,x_i,x_j,du,dv)\bigg)\bigg|^2\Bigg]
\\
&\leq \!C \|\psi\|_{H_x^{\e}\h^{\ktwo,\ttwo}}^4\!
\frac{1}{N^2M^2}\!\!\sum_{i,j=1}^N\!\sum_{k=1}^M\!\sup_{r\in[0,T]}\!\!\!\E\bigg[\!1\!+\!|\Bar{u}_{ik}^\epsilon|^{4\ttwo+4}\!+\!|\Bar{u}_{jk}^\epsilon|^{4\ttwo+4}\!+\!\int_{\R^{2\dirnum}}\!\!\!\!\!\!\!\!\!|u|^{4\ttwo+2}\!+\!|v|^{4\ttwo+4} df^{2,\epsilon}(r,x_i,x_j)\bigg]
\\
&\leq \frac{C}{M} \|\psi\|_{H_x^{\e}\h^{\ktwo,\ttwo}}^4
\bigg(1+\sup_{x\in Q}\E\big[|u_k(x,0)|^{4\ttwo+4}\big]\bigg).
\end{split}
\end{align}
}
In particular, we note that the last line vanishes in the $M,N\to\infty$ limit.
In the first passage we used convexity inequalities.
In the second passage we used that for fixed indices $i,j=1,\dots,N$ the couples $\{(\Bar{u}_{ik}^\epsilon,\Bar{u}_{jk}^\epsilon)\}_{k=1,\dots,M}$ are i.i.d. with law $f^{2,\epsilon}(x_i,x_j,t,du,dv)$, so that only diagonal terms survive when expanding the square of the sum over the index $k$.
In the third passage we used the definition of the operator $V_{(x,u)}^1$ and the estimate \eqref{formula/ norm of evaluation operators} on its norms, the linear growth properties \eqref{formula/ properties of b and sigma} of $\sigma$ and several applications of H\"older's and Young's inequality.
The last line follows from the moment estimates \eqref{formula/ uniform moment estimates}.

We finally consider the fifth term in \eqref{pointwise convergence of the quad variation/ proof 2}, which is \emph{deterministic}.
We start with the splitting:
{\small
\begin{align}\label{pointwise convergence of the quad variation/ proof 7}
\begin{split}
    E_r^5
    \!=\!&
    \!\!\sum_{i,j=1}^N\!\int_{Q_i^N\times Q_j^N}\!\! 
    \bigg( R^{\epsilon}(x_i-x_j)\int_{\R^{2\dirnum}}\nabla_u\psi(x_i,u)\sigma(x_i,r,u,f)\nabla_u\psi(x_j,v)\sigma(x_j,r,v,f) f^{2,\epsilon}(r,x_i,x_j,du,dv)
    \\
    &\qquad\qquad\qquad
    -
    R^{\epsilon}(x-y)\int_{\R^{2\dirnum}}\!\!\!\!\!\!\nabla_u\psi(x,u)\sigma(x,r,u,f)\nabla_u\psi(y,v)\sigma(y,r,v,f) f^{2,\epsilon}(r,x,y,du,dv)\bigg)\,dx\,dy
    \\
    =&
    \sum_{i,j=1}^N\!\int_{Q_i^N\times Q_j^N} \!\!\!
     \left(R^{\epsilon}(x_i-x_j)-R^{\epsilon}(x-y)\right)
     \\
     &\qquad\qquad\qquad\qquad
     \cdot\int_{\R^{2\dirnum}}\!\!\!\!\!\!\nabla_u\psi(x_i,u)\sigma(x_i,r,u,f)\nabla_u\psi(x_j,v)\sigma(x_j,r,v,f) f^{2,\epsilon}(r,x_i,x_j,du,dv)\,dx\,dy
     \\
     &+
     \sum_{i,j=1}^N\int_{Q_i^N\times Q_j^N} 
     \!\!\!\!\!\!\!\!\!R^{\epsilon}(x-y)\bigg(\int_{\R^{2\dirnum}\times\R^{2\dirnum}}\!\!\!\!\!\!\!\!\nabla_u\psi(x_i,u)\sigma(x_i,r,u,f)\nabla_u\psi(x_j,v)\sigma(x_j,r,v,f) 
     \\
    &\qquad\qquad\qquad-\nabla_u\psi(x,u)\sigma(x,r,u,f)\nabla_u\psi(y,v)\sigma(y,r,v,f)\,\,\pi_r(x_i,x_j,x,y,du,dv,du',dv')\bigg)\,dx\,dy
    \\
    &=
    E_r^{5.1}+E_r^{5.2},    
\end{split}
\end{align}
}
where $\pi_r(x_i,x_j,x,y,du,dv,du',dv')$ denotes an optimal pairing in $\w_2(\R^{2\dirnum})$ between the measures $f^{2,\epsilon}(r,x_i,x_j,du,dv)$ and $f^{2,\epsilon}(r,x,y,du,dv)$.

For the first term in \eqref{pointwise convergence of the quad variation/ proof 7} we compute, for a constant $C=C(T,b,\sigma,\rho)$,
{\small
\begin{align}\label{pointwise convergence of the quad variation/ proof 8}
\begin{split}
    \left|E_r^{5.1}\right|
    &\leq C
    \sum_{i,j=1}^N\!\int_{Q_i^N\times Q_j^N} \!\!\!
     \left|R^{\epsilon}(x_i-x_j)-R^{\epsilon}(x-y)\right|
     \\
     &\qquad\qquad\qquad\qquad
     \cdot\int_{\R^{2\dirnum}}\!\!\!\|\psi\|_{H_x^{\e}\h^{\ktwo,\ttwo}}^2\big(1+|u|^{2\ttwo+2}+|v|^{2\ttwo+2}\big)\,f^{2,\epsilon}(r,x_i,x_j,du,dv)\,dx\,dy
     \\
     &\leq C
     \|\psi\|_{H_x^{\e}\h^{\ktwo,\ttwo}}^2
     \bigg(1+\sup_{x\in Q}\E\big[|u_k(x,0)|^{2\ttwo+2}\big]\bigg)\,N^{-\nicefrac{\alpha}{d}}\epsilon^{-1}.
\end{split}
\end{align}
}
In the first inequality we used the estimate \eqref{formula/ norm of evaluation operators} on $V_{(x,u)}^1$, the linear growth \eqref{formula/ properties of b and sigma} of $\sigma$ and several applications of H\"older's and Young's inequality.
In the second passage we used the moment estimates \eqref{formula/ uniform moment estimates}, the Lipschitzianity of $R^{\epsilon}$, namely
\begin{equation}
    R^{\epsilon}(x)-R^{\epsilon}(y)= C_{\rho}\int_{\R^d}\big(\rho\big(\nicefrac{x}{\epsilon}+z\big)-\rho\big(\nicefrac{y}{\epsilon}+z\big)\big)\rho(z)\,dz\,\lesssim\min\left\{\frac{x-y}{\epsilon},1\right\},
\end{equation}
and the fact that $\text{diam}(Q_i^N)\simeq N^{\nicefrac{1}{d}}$ and $\meas(Q_i^N)=\frac{1}{N}$.

For the second term in \eqref{pointwise convergence of the quad variation/ proof 7} we compute, for a constant $C=C(T,b,\sigma)$,
{\small
\begin{align}\label{pointwise convergence of the quad variation/ proof 9}
\begin{split}
    \left|E_r^{5.2}\right|\!
    \leq\!\!
    \!\!\sum_{i,j=1}^N\!\int_{Q_i^N\times Q_j^N}\!\!\!\!\!\!\!\!\!\!\!\!\!\!\!\!\!dx\,dy\int_{\R^{4\dirnum}} \!\!\!\!\!\!\!\!\!\!\Scale[0.9]{d\pi_r(x_i,x_j,x,y)}
    \bigg(\!
    &\Scale[1]{|\nabla_u\psi(x_i,u)\!-\!\nabla_u\psi(x,u')||\sigma(x_i,u,r,f)\nabla_u\psi(x_j,v)\sigma(x_j,v,r,f)|}
    \\
    &+
    |\nabla_u\psi(x,u')||\sigma(x_i,u,r,f)-\sigma(x,u',r,f)||\nabla_u\psi(x_j,v)\sigma(x_j,v,r,f)|
    \\
    &+
    |\nabla_u\psi(x,u')\sigma(x,u',r,f)||\nabla_u\psi(x_j,v)-\nabla_u\psi(y,v')||\sigma(x_j,v,r,f)|
    \\
    &+
    |\nabla_u\psi(x,u')\sigma(x,u',r,f)\nabla_u\psi(y,v)||\sigma(x_j,v,r,f)-\sigma(y,v',r,f)|
    \bigg)
    \\
    \leq C
    \!\!\sum_{i,j=1}^N\!\int_{Q_i^N\times Q_j^N}\!\!\!\!\!\!\!\!\!\!\!\!\!\!\!dx\,dy\int_{\R^{4\dirnum}}
    \!\!\!\!\!\!\big(\|\psi\|^2_{H_x^{\e}\h^{\ktwo,\ttwo}}&+1\big)\big(1+|u|+|v|+|u'|+|v'|\big)^{2\ttwo+2}
    \\
    \cdot\big(|x_i-&x|^{\alpha}+|x_j-y|^{\alpha}+|u-u'|+|v-v'|\big)\,\pi_r(x_i,x_j,x,y,du,dv,du',dv')
    \\
    \\
    \leq C
    \big(1+\|\psi\|^2_{H_x^{\e}\h^{\ktwo,\ttwo}}\big)
    \bigg(1+\sup_{x\in Q}\E&\big[|u_k(x,0)|^{4\ttwo+4}\big]\bigg)^{\nicefrac{1}{2}}
    \\
    \sum_{i,j=1}^N\!\int_{Q_i^N\times Q_j^N}\!\!\!\!
    |x_i-x|^{\alpha}&+|x_j-y|^{\alpha}+\w_2(\R^{2\dirnum})(f^{2,\epsilon}(r,x_i,x_j),f^{2,\epsilon}(r,x,y))\,\,dx\,dy  
    \\
    \leq C
    \big(1+\|\psi\|^2_{H_x^{\e}\h^{\ktwo,\ttwo}}\big)
    \bigg(1+\sup_{x\in Q}\E&\big[|u_k(x,0)|^{4\ttwo+4}\big]\bigg)\big(N^{-\nicefrac{\alpha}{d}}+N^{-\nicefrac{1}{d}}\epsilon^{-1}\big). 
\end{split}
\end{align}
}
In the first line we simply used the boundedness of $R^{\epsilon}$ and added and subtracted several mixed terms.
In the second passage we used the Lipschitz and linear growth properties \eqref{formula/ properties of b and sigma} of $\sigma$, the operators $V_{(x,u)}^1$ and $V_{(x,u),(y,v)}^{1,\text{dif}}$ and the estimates \eqref{formula/ norm of evaluation operators} on their norms, and several applications of H\"older's and Young's inequality.
In the third passage we used H\"older's inequality, the fact that $\pi$ is an optimal pairing and the moment estimates \eqref{formula/ uniform moment estimates}.
In the last passage we used Proposition \ref{proposition/ holderianity of the single and joint empirical measure} and that $\text{diam}(Q_i^N)\simeq N^{\nicefrac{1}{d}}$ and $\meas(Q_i^N)=\frac{1}{N}$.

Combining \eqref{pointwise convergence of the quad variation/ proof 7}, \eqref{pointwise convergence of the quad variation/ proof 8} and \eqref{pointwise convergence of the quad variation/ proof 9}, we have that $\displaystyle\sup_{r\in[0,T]}|E_r^5|\to0$ as $M,N\to\infty$.
This concludes the proof.
\end{proof}

\begin{proof}[\textbf{Proof of Proposition \ref{proposition/ H valued martingale MtMN}}]
\noindent
Consider a Hilbert basis $\{\psi_p\}_{p\geq1}$ of $H_x^{\e}\h^{\kone,\tone}$, we compute, for a constant $C=C(T,b,\sigma)$,
\begin{align}\label{H valued martingale MtMN/ proof 1}
    \begin{split}
        &\sup_{M,N}\,\E \bigg[\sup_{t\in[0,T]}\|M_t^{MN}\|_{H_x^{-\e}\h^{-\kone,\tone}}^2\bigg]
        \\
        &\leq
        \sup_{M,N}\E\bigg[\sum_{p\geq1} \sup_{t\in[0,T]}|M_t^{MN}(\psi_p)|^2\bigg]
        \\
        &\leq C
        \sup_{M,N}\E\bigg[\sum_{p\geq1} \langle M^{MN}_T(\psi_p)\rangle\bigg]
        \\
        &\leq C
        \sup_{M,N}\frac{C_{MN}^2}{M}\frac{1}{MN^2}\sum_{k=1}^M\sum_{i,j=1}^N\int_0^T\E\bigg[\sum_{p\geq1}|\nabla_u\psi_p(x_i,u^\epsilon_{ik})|^2|\sigma(x_i,u_{ik}^\epsilon,r,f_{MN})|^2\bigg]\,dr\, R^\epsilon(x_i-x_j)
        \\
        &\leq C
        \sup_{M,N}\frac{C_{MN}^2}{M}\frac{1}{MN}\sum_{k=1}^M\sum_{i=1}^N\int_0^T1+\E\big[|u_{ik}^\epsilon|^{2\tone+2}\big]\,dr
        \\
        &\leq C
        \sup_{M,N}\frac{C_{MN}^2}{M}\bigg(1+\sup_{x\in Q}\E\big[|u_k(x,0)|^{2\tone+2}\big]\bigg).
    \end{split}
\end{align}
In the first passage we simply used Parseval identity and in the second passage we used the Burkholder-Davis-Gundy inequality.
In the third passage we used the definition \eqref{formula/ quadratic variation of real martingale term} of the quadratic variation and the Cauchy–Schwarz inequality.
In the fourth passage we used the Parseval identity $$\|V^1_{x_i,u^\epsilon_{ik}}\|_{H_x^{-\e}\h^{-\kone,\tone}}=\sum_{p\geq1}|\nabla_u\psi_p(x_i,u^\epsilon_{ik})|^2$$ together with the estimates \eqref{formula/ norm of evaluation operators}, the sublinear growth \eqref{formula/ properties of b and sigma} of $\sigma$ and the boundedness of $R^\epsilon$.
In the last passage we used the uniform moment estimates \eqref{formula/ uniform moment estimates}.

We now show the continuity of the trajectories.
By \eqref{H valued martingale MtMN/ proof 1} and Chebyshev's inequality we find a subset of full probability $\Omega_0\subset \Omega$ such that
\begin{equation}\label{H valued martingale MtMN/ proof 2}
\forall\omega\in\Omega_0\,\,\forall M,N\in\N\,\,\forall \delta>0\,\,\exists K=K(M,N,\omega,\delta)\,\,\text{s.t.}\,\sum_{p\geq K}\sup_{t\in[0,T]}|M_t^{MN}(\psi_p)(\omega)|^2<\delta,
\end{equation}
and such that $r\mapsto M_r^{MN}(\psi_p)(\omega)$ is continuous for each $p\in\N$.
Now for any $\omega\in\Omega_0$ and any $s,t\in[0,T]$, we compute
\begin{align}\label{H valued martingale MtMN/ proof 3}
    \begin{split}
        \|M_t^{MN}-M_s^{MN}\|_{H_x^{-\e}\h^{-\kone,\tone}}
        &=
        \sum_{p\geq1}|M_t^{MN}(\psi_p)-M_s^{MN}(\psi_p)|^2
        \\
        &\leq
        \sum_{p=1}^{K-1}|M_t^{MN}(\psi_p)-M_s^{MN}(\psi_p)|^2
        +
        2\sum_{p\geq K}\sup_{r\in[0,T]}|M_r^{MN}(\psi_p)|^2
        \\
        &\leq
        \sum_{p=1}^{K-1}|M_t^{MN}(\psi_p)-M_s^{MN}(\psi_p)|^2+2\delta.
    \end{split}
\end{align}
Continuity follows since $r\mapsto M_r^{MN}(\psi_p)$ is continuous for each $p\in\N$ and $\omega\in\Omega_0$.
\end{proof}

\begin{proof}[\textbf{Proof of Proposition \ref{proposition/ first bounds on etaMN}}]
\noindent
We start with the splitting
\begin{align}\label{first bounds on etaMN/ proof 1}
    \begin{split}
        \eta_t^{MN}=\cmn(f_{MN}-\Bar{f}_{MN})+\cmn(\Bar{f}_{MN}-\bar{f}_N)+\cmn(\bar{f}_N-f),
    \end{split}
\end{align}
for the measures
\begin{align}\label{first bounds on etaMN/ proof 2}
    \begin{split}
        \bar{f}_{MN}(t,dx,du)\coloneqq\frac{1}{MN}\sum_{j=1}^N\sum_{k=1}^M\delta_{(x_j,\Bar{u}_{jk}^\epsilon)}\,\text{ and }\,
        \bar{f}_N(t,dx,du)\coloneqq\frac{1}{N}\sum_{j=1}^N\delta_{x_j}\otimes f(t,x_j,du).
    \end{split}
\end{align}
For a Hilbert basis of $\{\psi_p\}_{p\geq1}$ of $H_x^{\e}\h^{\kone,\tone}$, Parseval identity and convexity inequalities yield, for a numeric constant $C$,
{\small
\begin{align}\label{first bounds on etaMN/ proof 3}
    \begin{split}
        \|\eta_t^{MN}\|_{H_x^{\e}\h^{\kone,\tone}}\leq C\sum_{p\geq1}\,\langle\cmn(f_{MN}-\Bar{f}_{MN}),\psi_p\rangle^2+
        \langle\cmn(\Bar{f}_{MN}-\bar{f}_N),\psi_p\rangle^2
        +
        \langle\cmn(\bar{f}_N-f),\psi_p\rangle^2.
    \end{split}
\end{align}}

We handle each sum in \eqref{first bounds on etaMN/ proof 3} separately.
The result will follow from \eqref{first bounds on etaMN/ proof 3} and estimates \eqref{first bounds on etaMN/ proof 4}, \eqref{first bounds on etaMN/ proof 5} and \eqref{first bounds on etaMN/ proof 6} below.
For the first one we compute, for a constant $C=C(T,b,\sigma)$,
{\small
\begin{align}\label{first bounds on etaMN/ proof 4}
    \begin{split}
        \E\bigg[\sum_{p\geq1}\,\langle\cmn(f_{MN}-\Bar{f}_{MN}),\psi_p\rangle^2\bigg]
        &
        \leq C\,\frac{\cmn^2}{MN}\sum_{i=1}^N\sum_{k=1}^M\E\bigg[\sum_{p\geq1}|\psi_p(x_i,u_{ik}^\epsilon)-\psi_p(x_i,\Bar{u}_{ik}^\epsilon)|^2\bigg]
        \\
        &\leq C\,\frac{\cmn^2}{MN}\sum_{i=1}^N\sum_{k=1}^M\E\bigg[\|V^{0,\text{dif}}_{(x_i,u_{ik}^\epsilon),(x_i,\Bar{u}_{ik}^\epsilon)}\|_{H_x^{-\e}\h^{-\kone,\tone}}^2\bigg]
        \\
        &\leq C\,\frac{\cmn^2}{MN}\sum_{i=1}^N\sum_{k=1}^M\E\big[|u_{ik}^\epsilon-\Bar{u}_{ik}^\epsilon|^4\big]^{\nicefrac{1}{2}}
        \E\big[1+|u_{ik}^\epsilon|^{4\tone}+|\Bar{u}_{ik}^\epsilon|^{4\tone}\big]^{\nicefrac{1}{2}}
        \\
        &\leq C\,\cmn^2 \left(\frac{1}{\sqrt{M}}+\frac{1}{N^{\nicefrac{\alpha}{d}}}\right)^2
        \bigg(1+\sup_{x\in Q}\E\big[|u_k(x,0)|^{4\tone}\big]^{\nicefrac{1}{2}}\bigg).
    \end{split}
\end{align}
}
In the first passage we used the definition of $f_{MN}$ and $\Bar{f}_{MN}$ and convexity properties.
In the second passage we used the definition \eqref{formula/ evaluation operators} of $V^{0,\text{dif}}_{(x_i,u_{ik}^\epsilon),(x_i,\Bar{u}_{ik}^\epsilon)}$ and Parseval identity.
In the third passage we exploited the estimate \eqref{formula/ norm of evaluation operators} on its norm and used the Cauchy–Schwarz inequality.
Finally in the last passage we used the error estimates \eqref{formula/ mean squared error estimates for actual particles vs McKean--Vlasov particles} and the moment estimates \eqref{formula/ uniform moment estimates}.

We now consider the second term in \eqref{first bounds on etaMN/ proof 3}.
For a constant $C=C(T,b,\sigma)$, we compute
{\small
\begin{align}\label{first bounds on etaMN/ proof 5}
    \begin{split}
        \E\bigg[\sum_{p\geq1}\,\langle\cmn(\bar{f}_{MN}-\Bar{f}_{N}),\psi_p\rangle^2\bigg]
        &
        \leq C\,\frac{\cmn}{N}\sum_{i=1}^N\frac{1}{M^2}\E\bigg[\sum_{p\geq1}\bigg|\sum_{k=1}^M\psi_p(x_i,\bar{u}_{ik}^\epsilon)-\int_{\R^{\dirnum}}\psi_p(x_i,v)f(x_i,t,dv)\bigg|^2\bigg]
        \\
        &
        \leq C\,\frac{\cmn^2}{N}\sum_{i=1}^N\frac{1}{M}\E\bigg[\sum_{p\geq1}\bigg|\psi_p(x_i,\bar{u}_{ik}^\epsilon)-\int_{\R^{\dirnum}}\psi_p(x_i,v)f(x_i,t,dv)\bigg|^2\bigg]
        \\
        &
        \leq C\,\frac{\cmn^2}{M}\frac{1}{N}\sum_{i=1}^N\E\big[\sum_{p\geq1}\big|\psi_p(x_i,\bar{u}_{ik}^\epsilon)\big|^2\big]
        \\
        &
        \leq C\,\frac{\cmn^2}{M}\frac{1}{N}\sum_{i=1}^N\E\big[\|V^{0}_{x_i,\bar{u}_{ik}^\epsilon}\|_{H_x^{-\e}\h^{-\kone,\tone}}^2\big]
        \\
        &
        \leq C\,\frac{\cmn^2}{M}      
        \bigg(1+\sup_{x\in Q}\E\big[|u_k(x,0)|^{2\tone}\big]\bigg).       
    \end{split}
\end{align}}
In the first passage we used the definition of $\Bar{f}_{MN}$ and $\Bar{f}_N$ and convexity.
In the second passage we used that for each fixed $i=1,\dots,N$ the particles $\{\Bar{u}_{ik}^\epsilon\}_{k=1,\dots,M}$ are i.i.d. with law $f(x_i,t,dv)$, so that only diagonal terms survive when expanding the square of the sum over the index $k$.
In the third passage we used $(a-b)^2\leq 2a^2+2b^2$, H\"older's inequality and again that the $\Bar{u}_{ik}^\epsilon$ are i.i.d. with law $f(x_i,t,dv)$.
In the fourth passage we used the definition \eqref{formula/ evaluation operators} of $V^{0}_{x_i,\Bar{u}_{ik}^\epsilon}$ and Parseval identity.
In the fifth passage we exploited the estimate \eqref{formula/ norm of evaluation operators} on its norm and the moment estimates \eqref{formula/ uniform moment estimates}.

Finally we consider the third term in \eqref{first bounds on etaMN/ proof 3}.
For every $x,y\in Q$, let $\pi_t(x_i,y,du,dv)$ denote an optimal pairing between $f(x,t,u)$ and $f(y,t,v)$ for the Wasserstein distance $\w_4(\R^{\dirnum})$.
For this \emph{deterministic} term we compute, for a constant $C=C(T,b,\sigma)$,
{\small
\begin{align}\label{first bounds on etaMN/ proof 6}
    \begin{split}
        \sum_{p\geq1}&\,\langle\cmn(\bar{f}_{N}-f),\psi_p\rangle^2
        \\
        &=
        \sum_{p\geq1}\cmn^2\bigg(\sum_{i=1}^N\int_{Q_i^N}\int_{\R^{2\dirnum}}\psi_p(x_i,u)-\psi_p(y,v)\,\pi_t(x_i,y,du,dv)\,dy\bigg)^2
        \\
        &\leq C
        \cmn^2\sum_{i=1}^N\int_{Q_i^N}\int_{\R^{2\dirnum}}\sum_{p\geq1}\left|\psi_p(x_i,u)-\psi_p(y,v)\right|^2\pi_t(x_i,y,du,dv)\,dy
        \\
        &\leq C
        \cmn^2\sum_{i=1}^N\int_{Q_i^N}\int_{\R^{2\dirnum}}  (1+|u|+|v|)^{2\tone}(|x_i-y|^{2\alpha}+|u-v|^2)\,\pi_t(x_i,y,du,dv)\,dy
        \\
        &\leq C
        \cmn^2\sum_{i=1}^N\int_{Q_i^N}\bigg(\int_{\R^{\dirnum}}  1+|u|^{4\tone}+|v|^{4\tone}\pi_t(x_i,y,du,dv)\bigg)^{\nicefrac{1}{2}}
        \\
        &\qquad\qquad\qquad\qquad\qquad\cdot
        \Big( |x_i-y|^{4\alpha}+\w_4(\R^{\dirnum})(f(x_i,t,du),f(y,t,dv)) \Big)^{\nicefrac{1}{2}}dy
        \\
        &\leq 
        C\cmn^2\Big(1+\sup_{x\in Q}\E\big[|u_k(x,0)|^{4\tone}\big]^{\nicefrac{1}{2}}\Big)N^{-\nicefrac{2\alpha}{d}}.    
    \end{split}
\end{align}
}
In the first passage we simply used the definitions of $\Bar{f}_N$, $f$ and $\pi_t(x_i,y)$ and that $\meas(Q_i^N)=\nicefrac{1}{N}$.
In the second passage we used Parseval identity for the quantity $\|V^{0,\text{dif}}_{(x_i,u),(y,v)}\|_{H_x^{-\e}\h^{-\kone,\tone}}^2$ and estimate \eqref{formula/ norm of evaluation operators}.
In the third inequality we used H\"older's inequality and the choice of $\pi_t$.
In the last passage we used Proposition \ref{proposition/ holderianity of the single and joint empirical measure}, the moment estimates \eqref{formula/ uniform moment estimates} and the fact that $\text{diam}(Q_i^N)\simeq N^{-\nicefrac{1}{d}}$.
\end{proof}

\begin{proof}[\textbf{Proof of Proposition \ref{proposition/ norm bounds on the linearized operator}}]
\noindent
We prove the first bound in \eqref{formula/ norm of linearized operators 1}, the other bounds are proved almost identically.
To streamline the proof, we consider two prototypical terms in the expression \eqref{formula/ linearized operator} for the operator $\lin(f_{MN},f)$, namely the second and the fourth.
The remaining terms are treated analogously, and are simpler since they involve less derivatives in $(x,u)$ or $(y,v)$ and less `unbounded products' of terms involving $b$, $\sigma$ and $\psi$.

For the second term in \eqref{formula/ linearized operator}, using the linear growth properties \eqref{formula/ properties of b and sigma} of the coefficient $\sigma$ and its derivative, we compute, for a constant $C=C(T,\sigma,d,\dirnum)$,
\begin{align}\label{norm bounds on the linearized operator/ proof 1}
    \begin{split}
    \big\|\Delta_v&\psi(y,v)\,\sigma(y,v,t,f_{MN})^2\big\|_{H_x^{\e}\h^{\kone,\tone}}^2
    \\
    &:=\sum_{j=0}^{\e}\sum_{h=0}^{\kone}\int_{Q}\int_{\R^{\dirnum}}\big|D_y^jD_v^h\big(\Delta_v\,\psi\sigma^2\big)\big|^2(1+|v|)^{-2\tone}dy\,dv
        \\
        &\leq C
        \int_{Q}\int_{\R^{\dirnum}}
        \bigg(\sum_{j=0}^{\e}\sum_{h=0}^{\kone}\left|D_y^jD_v^{h+2}\psi\right|^2\bigg) \bigg(\sum_{j=0}^{\e}\sum_{h=0}^{\kone}\left|D_y^jD_v^{h}\sigma^2\right|^2\bigg)(1+|v|)^{-2\tone}dy\,dv
        \\
        &\leq C \sum_{j=0}^{\e}\sum_{h=0}^{\kone+2}\int_{Q}\int_{\R^{\dirnum}}\left|D_y^jD_v^{h}\psi\right|^2(1+|v|)^{-2(\tone-2)}dv\,dy
        \\
        &=:C\|\psi\|_{H_x^{\e}\h^{\kone+2,\tone-2}}\leq C
        \|\psi\|_{H_x^{\e}\h^{\ktwo,\ttwo}}.
    \end{split}
\end{align}

Similarly for the fourth term in \eqref{formula/ linearized operator}, we compute, for a constant $C=C(T,\sigma,\dirnum,d)$,
{\small
\begin{align}\label{norm bounds on the linearized operator/ proof 2}
    \begin{split}
    \bigg\|&\int_{Q\times\R^{\dirnum}}\!\!\!\!\!\!\!\!\!\!\!\!\!\!\Delta_u\psi(x,u)\sigma_1(x,y,t,u,v)2\sigma_0(x,t,u)\Big(\text{\scalebox{0.9}{$\int_0^1\Dot{\phi}\big((1-\lambda)\sigma_1(x,t,u,f)+\lambda \sigma_1(x,t,u,f_{MN})\big)d\lambda$}}\Big)f(t,dx,du)\bigg\|_{H_x^{\e}\h^{\kone,\tone}}^2
    \\
    &\leq C\sum_{j=0}^{\e}\sum_{h=0}^{\kone}\int_{Q\times\R^{\dirnum}}\!\!\int_{Q\times\R^{\dirnum}}\!\!\!\!\!\!\!\!\!\!\!|\Delta_u\psi(x,u)|^2|D_y^jD_v^h\sigma_1(x,y,t,u,v)|^2|\sigma_0(x,t,u)|^2|\dot{\phi}|^2_{\infty}(1+|v|)^{-2\tone}f(t,dx,du)\,dy\,dv
    \\
    &\leq C 
    \int_{Q\times\R^{\dirnum}}\int_{Q\times\R^{\dirnum}}\|\psi\|_{H_x^{\e}\h^{\ktwo,\ttwo}}^2(1+|u|)^{4+2\ttwo}(1+|v|)^{2-2\tone}f(t,dx,du)\,dy\,dv
    \\
    &\leq C \,\|\psi\|_{H_x^{\e}\h^{\ktwo,\ttwo}}^2 \Big(1+\sup_{x\in Q}\E\big[|u_k(x,0)|^{4+2\ttwo}\big]\Big).
    \end{split}
\end{align}
}
In the first passage we simply used the definition of the norm $\|\cdot\|_{H_x^{\e}\h^{\kone,\tone}}$, differentiated under the integral sign and applied H\"older's inequality.
In the second passage we used the boundedness of $\dot{\phi}$, the linear growth assumptions \eqref{formula/ assumption on b0 sigma0 first}-\eqref{formula/ assumption on b1 sigma1 second} on $\sigma_0$ and $\sigma_1$ and the embedding \eqref{formula/ continuous embedding of working spaces into bounded functions} to bound $\Delta_u\psi$.
In the last passage, the integral in the variable $v$ is finite since $2(\tone-1)>\dirnum$ and we applied the moment estimates \eqref{formula/ uniform moment estimates} in the $u$ variable.
\end{proof}

\begin{proof}[\textbf{Proof of Proposition \ref{proposition/ H valued semimartingale etaMN}}]
\noindent
We already know that {\small$\eta_t^{MN}$} takes values in {\small$H^{-\e}\h^{-\kone,\tone}\subseteq H_x^{-\e}\h^{-\ktwo,\ttwo}$} and that for every $\psi\in H_x^{\e}\h^{\kone,\tone}$ the process $\langle \eta_t^{MN},\psi\rangle$ is a real coninuous semimartingale with decomposition \eqref{formula/ real semimartingale rewriting of eta epsilon}.
We are left with showing the bound \eqref{formula/ second bounds on etaMN} and the continuity of the trajectories in $H_x^{-\e}\h^{-\ktwo,\ttwo}$.
The decomposition \eqref{formula/ H valued semimartingale rewriting of etaMN} at the distribution level then follows immediately from the decomposition \eqref{formula/ real semimartingale rewriting of eta epsilon} and from the estimate \eqref{formula/ norm of linearized operators 1}, which implies that the integral in \eqref{formula/ H valued semimartingale rewriting of etaMN} is well-defined.
For the bound \eqref{formula/ second bounds on etaMN}, we consider a Hilbert basis $\{\psi_p\}_{p\geq1}$ of $H_x^{\e}\h^{\ktwo,\ttwo}$.
For any $M,N\in\N$ we compute, for a constant $C=C(T,b,\sigma,\dirnum,d)$,
\small
\begin{align}\label{H valued semimartingale rewriting of etaMN/ proof 1}
    \begin{split}
    \E &\bigg[\sup_{t\in[0,T]}\|\eta_t^{MN}\|_{H_x^{-\e}\h^{-\ktwo,\ttwo}}^2\bigg]
    \\
    &\leq
    \E\bigg[\sum_{p\geq1} \sup_{t\in[0,T]}|\langle\eta_t^{MN},\psi_p\rangle|^2\bigg]
    \\
    &\leq C \bigg(
    \E\Big[\|\eta_0^{MN}\|_{H_x^{-\e}\h^{-\ktwo,\ttwo}}^2\Big]
    + \E\Big[\|M_T^{MN}\|_{H_x^{-\e}\h^{-\ktwo,\ttwo}}^2\Big]
    \\
    &\qquad\qquad\quad+
    \int_0^T\E\Big[\|\lin_r(f_{MN},f)^*[\eta^{MN}_r]\|_{H_x^{-\e}\h^{-\ktwo,\ttwo}}^2\Big]\,dr\bigg)
    \\
    &\leq C\bigg(
    \E\Big[\|\eta_0^{MN}\|_{H_x^{-\e}\h^{-\kone,\tone}}^2\Big]
    + \E\Big[\|M_T^{MN}\|_{H_x^{-\e}\h^{-\kone,\tone}}^2\Big]
    \\
    &\qquad\qquad+
    \!\!\sup_{r\in[0,T]}\!\!\underset{\omega}{\esssup}\|\lin_r(f_{MN},f)^*\|_{H_x^{-\e}\h^{-\kone,\tone},H_x^{-\e}\h^{-\ktwo,\ttwo}}\sup_{r\in[0,T]}\E\big[\|\eta^{MN}_r\|_{H_x^{-\e}\h^{-\kone,\tone}}^2\big]\,dr\bigg)
    \\
    &\leq C
    \sup_{M,N}\left(\frac{\cmn^2}{M}+\frac{\cmn^2}{N^{\nicefrac{2\alpha}{d}}}\right)\bigg(1+\sup_{x\in Q}\E\Big[|u(x,0)|^{4\tone}\Big]\bigg).
    \end{split}
\end{align}
In the first passage we used Parseval identity.
In the second passage we exploited the decomposition \eqref{formula/ real semimartingale rewriting of eta epsilon} for each $\psi_p$, convexity and H\"older's inequality, and again several applications of the Parseval identity.
In the last passage we simply inserted the estimates \eqref{formula/ uniform estimate on martingale term}, \eqref{formula/ first bounds on etaMN} and \eqref{formula/ norm of linearized operators 1}.

The continuity of the trajectories in $H_x^{\-e}\h^{-\ktwo,\ttwo}$ now follows from \eqref{H valued semimartingale rewriting of etaMN/ proof 1} using the exact same argument as in\eqref{H valued martingale MtMN/ proof 2}-\eqref{H valued martingale MtMN/ proof 3}, with $\eta_t^{MN}$ and $H_x^{\-e}\h^{-\ktwo,\ttwo}$ replacing $M_t^{MN}$ and $H_x^{\-e}\h^{-\kone,\tone}$.
\end{proof}

\begin{proof}[\textbf{Proof of Proposition \ref{proposition/ tightness of the processes}}]
\noindent
We first take care of the martingale term.
Thanks to the Ascoli--Arzelà criterion for compactness, we have to show that $M_t^{MN}$ satisfies the following conditions.
\begin{itemize}
    \item[(C1)] For every $t$ in a dense subset of $\R^+$ the laws of the r.v. $M_t^{MN}$ are tight in $H_x^{-\eplus}\h^{-\ktwo,\ttwo}$.
    
    \item[(C2)] For every $T\geq0$ and every $\beta,\lambda>0$, there exists $\delta>0$ such that 
        \begin{equation}\label{tightness of the processes/ proof 1}
        \limsup_{M,N\to\infty}\p\bigg(\sup_{\substack{t,s\in[0,T],\,|t-s|\leq\delta}}\|M_t^{MN}-M_s^{MN}\|_{H_x^{-\eplus}\h^{-\ktwo,\ttwo}}>\beta\bigg)<\lambda.
        \end{equation}
\end{itemize}

Condition (C1) is readily verified.
Indeed, the compact embedding $H_x^{-e}\h^{-\kone,\tone}\doublehookrightarrow H_x^{-\eplus\h^{-\ktwo,\ttwo}}$, which follows from \eqref{formula/ compact embedding between working spaces}, ensures that balls of $H_x^{-e}\h^{-\kone,\tone}$ are compact in $H_x^{-\eplus\h^{-\ktwo,\ttwo}}$.
In turn, the uniform estimate \eqref{formula/ uniform estimate on martingale term} and Chebyshev's inequality yield, for $C=C(T,b,\sigma)$,
\begin{equation}\label{tightness of the processes/ proof 2}
    \sup_{M,N}\p\Big(\|M_t^{MN}\|_{H_x^{-\e}\h^{-\kone,\tone}}\geq R\Big)
    \leq
    \frac{C}{R^2}\sup_{M,N}\E\bigg[\sup_{t\in[0,T]}\|M_t^{MN}\|_{H_x^{-\e}\h^{-\kone,\tone}}^2\bigg]
    \leq
    \frac{C}{R^2}\sup_{M,N}\frac{\cmn^2}{M}.
\end{equation}

Condition (C2) is verified using the well-known Aldous criterion (see e.g. \cite[Theorem 2.2.2]{Joffe_metivier_weak_convergence_of_sequence_of_semimartingales}).
In our context, it establishes that condition (C2) is implied by the following condition.
\begin{itemize}
   \item[(A)] For every $T\geq0$ and every $\beta,\lambda>0$, there exist $\delta>0$ and $K\in\N$ such that, for every sequence $\tau_{MN}$ of $\mathcal{F}_t$-stopping times satisfying $\tau_{MN}\leq T$ almost surely, we have
        \begin{equation}\label{tightness of the processes/ proof 3}
        \sup_{M,N\geq K}\sup_{\theta\leq \delta}
        \p\Big(\|M_{\tau_{MN}+\theta}^{MN}-M_{\tau_{MN}}^{MN}\|_{H_x^{-\eplus}\h^{-\ktwo,\ttwo}}>\beta\Big)<\lambda.
        \end{equation}
\end{itemize}
Condition (A) follows immediately from the next inequality.
Let $\{\psi_p\}_{p\geq 1}$ be a Hilbert basis of $H_x^{\e}\h^{\ktwo,\ttwo}$.
We compute, for a constant $C=C(T,b,\sigma,\dirnum,d)$,
{\small
\begin{align}\label{H valued martingale MtMN/ proof 4}
    \begin{split}
        \p&\Big(\|M_{\tau_{MN}+\theta}^{MN}-M_{\tau_{MN}}^{MN}\|_{H_x^{-\eplus}\h^{-\ktwo,\ttwo}}>\beta\Big)
        \\
        &\leq
        \frac{1}{\beta^2}\E\bigg[\sum_{p\geq1}\big|M_{\tau_{MN}+\theta}^{MN}(\psi_p)-M_{\tau_{MN}}^{MN}(\psi_p)\big|^2\bigg]
        \\
        &\leq
        \frac{C}{\beta^2}\frac{\cmn^2}{M}\E\bigg[\sum_{p\geq1}\int_{\tau_{MN}}^{\tau_{MN}+\theta}\!\!\!\int_{Q^2\times\R^{2\dirnum}}\!\!\!\!\!\!\big(|\nabla_u\psi_p(x,u)|^2+|\nabla_v\psi_p(y,v)|^2\big)(1+|u|)(1+|v|)\,f_{MN}^2(r,dx,dy,du,dv)\,dr\bigg]
        \\
        &\leq
        \frac{C}{\beta^2}\frac{\cmn^2}{M}\E\bigg[\delta \sup_{r\in[0,T]}\int_{Q^2\times\R^{2\dirnum}}\!\!\!\!\!\!1+|u|^{2\ttwo+2}+|v|^{2\ttwo+2}\,\,f_{MN}^2(r,dx,dy,du,dv)\bigg]
        \\
        &\leq 
        \frac{\delta}{\beta^2}C\sup_{M,N}\frac{\cmn^2}{M}\Big(1+\sup_{x\in Q}\E\big[|u_k(x,0)|^{2\ttwo+2}\big]\Big).
    \end{split}
\end{align}
}
In the first line we used Parseval identity and Chebyshev's inequality.
In the second line we used the Burkholder-Davis-Gundy inequality and the expression \eqref{formula/ quadratic variation of real martingale term} for the quadratic variation, and then the Cauchy–Schwarz inequality, the linear growth \eqref{formula/ properties of b and sigma} of $\sigma$ and the boundedness of $R^{\epsilon}$.
In the third passage we used Parseval identity for the norm $\|V_{x,u}^1\|_{H_x^{-\e}\h^{-\ktwo,\ttwo}}$ and the bound \eqref{formula/ norm of evaluation operators} and Young's and H\"older's inequality.
In the last passage we used the definition of $f_{MN}^2$ and the moment estimates \eqref{formula/ uniform moment estimates}.

We now prove the tightness of the fluctuations $\eta_t^{MN}$.
For this, we prove that each term of the semimartingale decomposition \eqref{formula/ H valued semimartingale rewriting of etaMN} is tight.
We have just proved that the term $M_t^{MN}$ is tight in $C([0,\infty);H_x^{-\eplus}\h^{-\ktwo,\ttwo}\big)$.
Furthermore, as in \eqref{H valued martingale MtMN/ proof 2}, condition (C1) holds for $\eta_t^{MN}$ thanks to the compact embedding $H_x^{-e}\h^{-\kone,\tone}\doublehookrightarrow H_x^{-\eplus\h^{-\ktwo,\ttwo}}$ and to estimate \eqref{formula/ first bounds on etaMN}.
In particular, condition (C1) holds for $\eta_0^{MN}$ and in turn also for the integral term in \eqref{formula/ H valued semimartingale rewriting of etaMN}.
Since the initial data is independent of time, condition (C2) is automatically satisfied and $\eta_0^{MN}$ is tight.
To conclude, we just have to show that the integral term in \eqref{formula/ H valued semimartingale rewriting of etaMN} satisfies condition (C2).
For this we compute, for a constant $C=C(T,b,\sigma)$,
\begin{align}\eqref{H valued martingale MtMN/ proof 5}
    \begin{split}
        \p&\bigg(\sup_{|s-t|\leq\delta}\Big\|\int_s^t\lin_r(f_{MN},f)^*[\eta^{MN}_r]\,dr\Big\|_{H_x^{-\eplus}\h^{-\ktwo,\ttwo}}>\beta\bigg)
        \\
        &\leq 
        \p\bigg(\sup_{|s-t|\leq\delta}\int_s^t\|\lin_r(f_{MN},f)\|_{H_x^{\e}\h^{\ktwo,\ttwo},H_x^{\e}\h^{\kone,\tone}}\|\eta^{MN}_r\|_{H_x^{-\e}\h^{-\kone,\tone}}dr>\beta\bigg)
        \\
        &\leq 
        \p\bigg(\delta^{\nicefrac{1}{2}}\sup_{t\in[0,T]}\underset{\omega\in\Omega}{\esssup}\|\lin_t(f_{MN},f)\|_{H_x^{\e}\h^{\ktwo,\ttwo},H_x^{\e}\h^{\kone,\tone}}\Big(\int_0^T\|\eta^{MN}_r\|^2_{H_x^{-\e}\h^{-\kone,\tone}}dr\Big)^{\nicefrac{1}{2}}>\beta\bigg)
        \\
        &\leq
        \frac{\delta}{\beta^2}\sup_{t\in[0,T]}\underset{\omega\in\Omega}{\esssup}\|\lin_t(f_{MN},f)\|_{H_x^{\e}\h^{\ktwo,\ttwo},H_x^{\e}\h^{\kone,\tone}}^2\E\Big[\int_0^T\|\eta^{MN}_r\|^2_{H_x^{-\e}\h^{-\kone,\tone}}dr\Big]
        \\
        &\leq
        \frac{\delta}{\beta^2}\sup_{MN}\bigg(\sup_{t\in[0,T]}\underset{\omega\in\Omega}{\esssup}\|\lin_t(f_{MN},f)\|_{H_x^{\e}\h^{\ktwo,\ttwo},H_x^{\e}\h^{\kone,\tone}}^2
        T\sup_{r\in[0,T]}\E\bigg[\|\eta_r^{MN}\|^2_{H^{-\e}\h^{-\kone,\tone}}\bigg]\bigg).
    \end{split}
\end{align}
Now condition (C2) follows by using the estimates \eqref{formula/ norm of linearized operators 1} and \eqref{formula/ first bounds on etaMN} and choosing $\delta$ suitably small.
\end{proof}

\begin{proof}[\textbf{Proof of Proposition \ref{proposition/ clt for the initial data}}]
\noindent
By Proposition \ref{proposition/ tightness of the processes} we know that $\eta_0^{MN}=\cmn(f_{MN}(0)-f(0))$ is tight in $H_x^{-\eplus}\h^{-\ktwo,\ttwo}$.
For any $\psi\in H_x^{\e}\h^{\ktwo,\ttwo}$, we consider the splitting
\begin{align}\label{clt for the initial data/ proof 1}
\begin{split}
    \langle \eta_0^{MN},\psi \rangle
    =&
    \cm\bigg(\frac{1}{MN}\sum_{i=1}^N\sum_{k=1}^M\psi(x_i,u_k(x_i,0))-\int_Q\int_{\R^{\dirnum}}\psi(x,v)\,f_0(x,dv)\bigg)
    \\
    =&
    \frac{\cm}{M}\sum_{k=1}^M\bigg(\frac{1}{N}\sum_{i=1}^N\psi(x_i,u_k(x_i,0))-\int_{\R^{\dirnum}}\psi(x_i,v)\,f_0(x_i,dv)\bigg)
    \\
    &+
    \cm\sum_{i=1}^N\int_{Q_i^N}\bigg(\int_{\R^{\dirnum}}\psi(x_i,v)\,f_0(x_i,dv)-\int_{\R^{\dirnum}}\psi(x,v)\,f_0(x,dv)\bigg)\,dx
    \\
    =& Y_1+Y_2.
\end{split}
\end{align}

We claim that $Y_2\to0$ as $M,N\to\infty$ along the scaling $\cm N^{-\nicefrac{\alpha}{d}}\to0$.
Indeed we compute, for $\pi(x_i,x,du,dv)$ an optimal pairing in $\w_2(\R^{\dirnum})$ between $f_0(x_i,du)$ and $f_0(x,dv)$, for a constant $C=C(f_0)$,
{\small
\begin{align}\label{clt for the initial data/ proof 2}
\begin{split}
     |Y_2|
     &\leq  
     \cm\sum_{i=1}^N\int_{Q_i^N}\int_{\R^{2\dirnum}}|\psi(x_i,u)-\psi(x,v)|\,\pi(x_i,x,du,dv)\,dx
     \\
     &\leq
     C\|\psi\|_{H_x^{\eplus}\h^{\ktwo,\ttwo}}
     \cm\sum_{i=1}^N\int_{Q_i^N}\int_{\R^{2\dirnum}}
    (1+|u|+|v|)^{\ttwo}(|x_i-x|^{\alpha}+|u-v|)
     \,\pi(x_i,x,du,dv)\,dx
     \\
     &\leq
     C\|\psi\|_{H_x^{\eplus}\h^{\ktwo,\ttwo}}\bigg(1+\sup_{x\in Q}\E\big[|u_k(x,0)|^{2\ttwo}\big]\bigg)^{\nicefrac{1}{2}}
     \!\cm\!\sum_{i=1}^N\int_{Q_i^N}\!\!|x_i-x|^{\alpha}+\w_2(\R^{\dirnum})(f_0(x_i),f_0(x))\,dx
     \\
     &\leq
     C\|\psi\|_{H_x^{\eplus}\h^{\ktwo,\ttwo}}\bigg(1+\sup_{x\in Q}\E\big[|u_k(x,0)|^{2\ttwo}\big]\bigg)^{\nicefrac{1}{2}}\cm N^{-\nicefrac{\alpha}{d}}.
\end{split}
\end{align}
}
In the second passage we used the estimate \eqref{formula/ norm of evaluation operators} on the operator $V^{0,\text{dif}}_{(x_i,u),(x,v)}$.
In the third passage we used H\"older's inequality, the definition of $\pi$ and the moment estimates \eqref{formula/ uniform moment estimates}.
In the last passage we used the H\"older continuity of $f_0$ and that $\text{diam}(Q_i^N)\simeq N^{-\nicefrac{1}{d}}$ and $\meas(Q_i^N)=\frac{1}{N}$.

We now consider the first term in \eqref{clt for the initial data/ proof 1}.
We claim that its characteristic function satisfies 
\begin{equation}\label{clt for the initial data/ proof 3}
    \lim_{M,N\to\infty}\E[\exp(itY_1)]=e^{-\frac{t^2}{2}\mathcal{Q}(\psi,\psi)}.
\end{equation}
Combining this with \eqref{clt for the initial data/ proof 1} and \eqref{clt for the initial data/ proof 2} implies that, along a scaling regime such that $\sqrt{M}N^{-\nicefrac{\alpha}{d}}\to0$,
\begin{equation}\label{clt for the initial data/ proof 4}
\langle\eta_0^{MN},\psi\rangle\to\mathcal{N}\big(0,\mathcal{Q}(\psi,\psi)\big)\quad\text{in law for every $\psi\in H_x^{\e}\h^{\ktwo,\ttwo}$}.    
\end{equation}
Hence, using the tightness of $\eta_0^{MN}$ in $H_x^{-\eplus}\h^{-\ktwo,\ttwo}$, we conclude that
\begin{equation}\label{clt for the initial data/ proof 5}
\eta_0^{MN}\to\mathcal{N}\big(0,\mathcal{Q}\big)\quad\text{in law in $H_x^{-\eplus}\h^{-\ktwo,\ttwo}$}.    
\end{equation}

We prove \eqref{clt for the initial data/ proof 3}.
Recalling that the collections $\{u_k(x_i,0)\}_{i=1,\dots,N}$ are i.i.d. for $k=1,\dots,M$, we compute
\begin{align}\label{clt for the initial data/ proof 5}
\begin{split}
    \E[\exp(itY_1)]
    &=\bigg(\E\bigg[\frac{it}{\sqrt{M}}\frac{1}{N}\sum_{i=1}^N\psi(x_i,u_k(x_i,0))-\int_{\R^{\dirnum}}\psi(x_i,v)\,f_0(x_i,dv)\bigg]\bigg)^M
    \\
    &=
    \exp\bigg(-\frac{t^2}{2}\,\E\Big[\Scale[0.95]{\Big(\frac{1}{N}\sum_{i=1}^N\psi(x_i,u_k(x_i,0))-\int_{\R^{\dirnum}}\psi(x_i,v)\,f_0(x_i,dv)\Big)^2}\Big]+O\big(M^{-\nicefrac{1}{2}}\big)\bigg),
\end{split}
\end{align}
where the expression $O\big(M^{-\nicefrac{1}{2}}\big)$ depends on $\psi$ and $(u_k(x,0))_{x\in Q}$, but is uniform in $N\in\N$ since 
{\small
\begin{align}\label{clt for the initial data/ proof 4}
\begin{split}
    \Big|\frac{1}{N}\sum_{i=1}^N\psi(x_i,u_k(x_i,0))-\int_{\R^{\dirnum}}\psi(x_i,v)\,f_0(x_i,dv)\Big|^p
    \leq
    \sup_{x\in Q}\Big|\psi(x,u_k(x_i,0))-\int_{\R^{\dirnum}}\psi(x,v)\,f_0(x,dv)\Big|^p \,\,\,\forall p\geq1.
\end{split}
\end{align}
}
Therefore formula \eqref{clt for the initial data/ proof 3} will follow if we prove that
\begin{align}\label{clt for the initial data/ proof 5}
\begin{split}
    \lim_{N\to\infty}\E&\bigg[\Scale[1.1]{\Big(\frac{1}{N}\sum_{i=1}^N\psi(x_i,u_k(x_i,0))-\int_{\R^{\dirnum}}\psi(x_i,v)\,f_0(x_i,dv)\Big)^2}\bigg]
    \\
    &=\E\bigg[\Scale[1.1]{\Big(\int_Q\psi(x,u_k(x,0))-\int_{\R^{\dirnum}}\psi(x,v)\,f_0(x,dv)\,dx\Big)^2}\bigg]
    =\mathcal{Q}(\psi,\psi).
\end{split}
\end{align}
This follows from the next inequality.
For $\pi(x_i,x,du,dv)$ an optimal pairing in $\w_2(\R^{\dirnum})$ between $f_0(x_i,du)$ and $f_0(x,dv)$, for a constant $C=C(f_0,d,\dirnum)$, we compute
{\small
\begin{align}\label{clt for the initial data/ proof 6}
\begin{split}
 \E&\bigg[\Big|\frac{1}{N}\sum_{i=1}^N\Big(\psi(x_i,u_k(x_i,0))-\int_{\R^{\dirnum}}\!\!\!\psi(x_i,v)\,f_0(x_i,dv)\Big)-\Big(\int_Q\psi(x,u_k(x,0))-\int_{\R^{\dirnum}}\!\!\!\psi(x,v)\,f_0(x,dv)\,dx\Big)\Big|^2\bigg]
 \\
 &\leq
 \sum_{i=1}^N\int_{Q_i^N}\!\!\!\!\!dx\,\,\E\bigg[\Big|\,V^{0,\text{dif}}_{(x_i,u_k(x_i,0)),(x,u_k(x,0))}(\psi)\,+\int_{\R^{2\dirnum}}V^{0,\text{dif}}_{(x_i,u_k(x_i,0)),(x,u_k(x,0))}(\psi)\,\,\,\pi(x_i,x,du,dv)\Big|^2\bigg]
 \\
 &\leq
 C\|\psi\|_{H_x^{\e}\h^{\ktwo,\ttwo}}\Big(1+\sup_{x\in Q}\E\big[|u_k(x,0)|^{2\ttwo}\big]^{\nicefrac{1}{2}}\Big)
 \\
 &\qquad\qquad\qquad
 \cdot\sum_{i=1}^N\int_{Q_i^N}\!|x_i-x|^{\alpha}+\E\big[|u_k(x_i,0)-u_k(x,0)|^2\big]^{\nicefrac{1}{2}}+\w_2(\R^{\dirnum})(f_0(x_i),f_0(x))\,\,dx
 \\
 &\leq
 C\|\psi\|_{H_x^{\e}\h^{\ktwo,\ttwo}}\Big(1+\sup_{x\in Q}\E\big[|u_k(x,0)|^{2\ttwo}\big]^{\nicefrac{1}{2}}\Big)\, N^{-\nicefrac{\alpha}{d}}.
\end{split}
\end{align}
}
The first inequality follows from convexity.
In the second passage we first used the estimate \eqref{formula/ norm of evaluation operators} on the operators $V^{0,\text{dif}}$ and then used H\"older's inequality, the definition of $\pi$ and the moment estimates \eqref{formula/ uniform moment estimates}.
In the last passage we used the H\"older continuity of $f_0\in C^{\alpha}(Q;\p_2(\R^{\dirnum})$ and $u_k(\cdot,0)\in C(Q;L^2(\Omega))$, and that $\text{diam}(Q_i^N)\simeq N^{-\nicefrac{1}{d}}$ and $\meas(Q_i^N)=\frac{1}{N}$.
\end{proof}

\begin{proof}[\textbf{Proof of Theorem \ref{theorem/ clt/ clt section}}]
\noindent
By Proposition \ref{proposition/ tightness of the processes} the fluctuations {\small$\eta_t^{MN}$} are tight in {\small$C([0,T];H_x^{-\eplus}\h^{-\ktwo,\ttwo})$}.
Consider any subsequence $(M_k,N_k)$, which we relabel $(M,N)$, such that $\eta_t^{MN}$ is converging in law to some limit element $\eta_t^{\infty}$.
By Theorem \ref{theorem/ convergence of the single and joint empirical measure} and Proposition \ref{proposition/ clt for the martingale} and \ref{proposition/ clt for the initial data}, the sequences $\eta_0^{MN}$, $\eta_t^{MN}$, $M_t^{MN}$, $f_{MN}$ and $f_{MN}^2$ are converging in law to $\eta_0$, $\eta_t^{\infty}$, $G_t^{\epsilon}$, $f$ and $f^{2,\epsilon}$ in the space
{\small
\begin{equation}
    \label{clt/ proof 1}
    \mathbf{X}=H_x^{-\eplus}\h^{-\ktwo,\ttwo}\times C\big([0,T]; H_x^{-\eplus}\h^{-\ktwo,\ttwo}\big)^2\times C\big([0,T];\pr_2(Q\times\R^{\dirnum})\big)\times C\big([0,T];\pr_2(Q^2\times\R^{2\dirnum})\big).
\end{equation}
}
Using Skorokhod's representation theorem, we find another probability space $\Tilde{\Omega}$ and random sequences and limit elements supported in there such that
\begin{equation}
    \label{clt/ proof 2}
(\tilde{\eta}_0^{MN},\tilde{\eta}_t^{MN},\tilde{M}_t^{MN},\tilde{f}_{MN},\tilde{f}_{MN}^2)\to(\tilde{\eta}_0,\tilde{\eta}_t^{\infty},\tilde{G}_t^{\epsilon},f,f^{2,\epsilon})\text{ in $\mathbf{X}$, almost surely in $\Tilde{\Omega}$},
\end{equation}
and such that each term has the same law as the corresponding term supported in $\Omega$.
In particular, $\Tilde{\eta}_0$ is a Gaussian r.v. with mean zero and covariance $\mathcal{Q}$ given in \eqref{formula/ clt covariance of the limit initial data}, $\Tilde{G}_t^\epsilon$ is a Gaussian process with mean zero and covariance function $g_t^\epsilon$ given in \eqref{formula/ quadratic variation of real gaussian limit} and $f$, $f^{2,\epsilon}$ are the previous \emph{deterministic} measures.

Fatou's lemma and estimate \eqref{formula/ second bounds on etaMN} yield $\E\Big[\sup_{r\in[0,T]}\|\tilde{\eta}_r^{\epsilon}\|^2_{H_x^{-\eplus}\h^{-\ktwo,\ttwo}}\Big]<\infty$.
In particular, thanks to \eqref{formula/ norm of linearized operators 1}, the integral $\int_0^t\lin_r(f,f)^*\Tilde{\eta}^\epsilon_r\,dr$ is well-defined.
Furthermore, we still have the equation, written in weak form for $\psi\in H^{\eplus}\h^{\ktwo+2,\ttwo-2}$, for every $t\in[0,T]$, almost surely in $\Tilde{\Omega}$,
\begin{equation}\label{clt/ proof 3}
    \langle \tilde{\eta}_t^{MN},\psi\rangle-\int_0^t\langle\tilde{\eta}_r^{MN},\lin_r(\tilde{f}_{MN},f)[\psi]\,\rangle\,\, dr=\langle \tilde{\eta}_0^{MN}, \psi\rangle + \langle \tilde{M}_t^{MN},\psi\rangle.
\end{equation}
We claim that, along sub-subsequences, equation \eqref{clt/ proof 3} passes to the limit $M,N\to\infty$ and converges to
\begin{equation}\label{clt/ proof 4}
    \tilde{\eta}_t^{\epsilon}-\int_0^t\lin_r(f,f)^*[\tilde{\eta}_r^{\epsilon}]\,\, dr= \tilde{\eta}_0^{\epsilon} + \Tilde{G}_t^{\epsilon}.
\end{equation}
That is, it converges to the Langevin SPDE \eqref{formula/ langevin spde / clt section} written in $\Tilde{\Omega}$.

Assume the claim \eqref{clt/ proof 3}-\eqref{clt/ proof 4}.
It follows that $\Tilde{\eta}^{\epsilon}_t$ must be the unique-in-law weak solution in $C([0,T];H^{-\eplus}_x\h^{-(\ktwo+2),\ttwo-2})$ of the Langevin SPDE \eqref{clt/ proof 4}.
Since by construction $\eta_t^{\infty}$ and $\Tilde{\eta}_t^{\epsilon}$ have the same law in $C([0,T];H^{-\eplus}_x\h^{-\ktwo,\ttwo})$ and since $\eta_t^{MN}\to\eta_t^{\infty}$ in law, we conclude that along the initial subsequence $\eta_t^{MN}$ is converging in law to the solution of the SPDE \eqref{formula/ langevin spde / clt section}.
Since the converging subsequence $(M_k,N_k)$ was arbitrary, we conclude that the whole sequence $\eta_t^{MN}$ is converging to the solution of \eqref{formula/ langevin spde / clt section} and the theorem is proved.

We now prove the claim \eqref{clt/ proof 3}-\eqref{clt/ proof 4}.
By \eqref{clt/ proof 2} we know that $\Tilde{\eta}_0^{MN}\to\Tilde{\eta}_0$ and $\Tilde{M}_t^{MN}\to\Tilde{G}_t^{\epsilon}$ almost surely.
Hence the right-hand side of \eqref{clt/ proof 3} does converge to the right-hand side of \eqref{clt/ proof 4}.

We consider the left-hand side of the equations \eqref{clt/ proof 3}-\eqref{clt/ proof 4}.
Thanks to estimate \eqref{formula/ norm of linearized operators 1}, for every $\psi\in H^{\eplus}_x\h^{\ktwo+2,\ttwo-2}$ and $t\in[0,T]$, we have a continuous map
\begin{equation}\label{clt/ proof 5}
    F_{\psi}^t:C([0,T];H^{-\eplus}_x\h^{-\ktwo,\ttwo})\to\R\,\mid\,F_{\psi}^t(\zeta)=\langle \zeta,\psi\rangle-\int_0^t\langle\zeta,\lin_r(f,f)[\psi]\,\rangle\,\, dr.
\end{equation}
Since we know that $\Tilde{\eta}_{\cdot}^{MN}\to\Tilde{\eta}^{\epsilon}_{\cdot}$ almost surely in $C([0,T];H^{-\eplus}_x\h^{-\ktwo,\ttwo})$ and $F_{\psi}^t$ is continuous, we obtain that $F_{\psi}^t(\Tilde{\eta}_{\cdot}^{MN})\to F_{\psi}^t(\Tilde{\eta}^{\epsilon}_{\cdot})$ almost surely.
Hence the claim is proved upon showing that, almost surely along sub-subsequences,
\begin{equation}\label{clt/ proof 6}
    \int_0^t\langle\Tilde{\eta}^{MN}_r,\lin_r(\Tilde{f}_{MN},f)[\psi]\,\rangle\,\, dr-\int_0^t\langle\Tilde{\eta}^{MN}_r,\lin_r(f,f)[\psi]\,\rangle\,\, dr\to0.
\end{equation}
For this, we compute
\begin{align}\label{clt/ proof 7}
\begin{split}
    \E\Bigg[&\bigg|\int_0^t\langle\Tilde{\eta}^{MN}_r,\lin_r(\Tilde{f}_{MN},f)[\psi]\,\rangle\,\, dr-\int_0^t\langle\Tilde{\eta}^{MN}_r,\lin_r(f,f)[\psi]\,\rangle\,\, dr\bigg|\Bigg]
    \\
    &\leq
    \E\bigg[\sup_{r\in[0,T]}\|\eta_r^{MN}\|^2_{H_x^{-\e}\h^{-\ktwo,\ttwo}}\bigg]^{\nicefrac{1}{2}}\int_0^t
    \E\Big[\big\|(\lin_r(f_{MN},f)-\lin_r(f,f))[\psi]\big\|^2_{H_x^{\e}\h^{\ktwo,\ttwo}}\Big]^{\nicefrac{1}{2}} dr,
\end{split}
\end{align}
where we can drop the $\Tilde{\cdot}$ symbols and go back to working in $\Omega$ because of the equality in law.
The first term on the right-hand side of \eqref{clt/ proof 7} is bounded uniformly in $M,N$ by the estimates \eqref{formula/ second bounds on etaMN} and we finally show that integrand in the second term vanishes uniformly in time as $M,N\to\infty$.

This fact is proved almost identically to Lemma \ref{proposition/ norm bounds on the linearized operator}.
We show how to handle one of the terms in $\big\|(\lin_r(\Tilde{f}_{MN},f)-\lin_r(f,f))[\psi]\big\|^2_{H_x^{\e}\h^{\ktwo,\ttwo}}$, the others are treated analogously.
Namely, for the difference of the second terms in the expression \eqref{formula/ linearized operator} for $\lin_r(f_{MN},f)$ and $\lin_r(f,f)$, using the linear growth and Lipschitz properties \eqref{formula/ properties of b and sigma}, we compute, for a constant $C=(T,b,\sigma,d,\dirnum)$,
{\small
\begin{align}\label{clt/ proof 8}
    \begin{split}
    \E&\Big[\big\|\Delta_v\psi(y,v)\,\sigma(y,v,t,f_{MN})^2-\Delta_v\psi(y,v)\,\sigma(y,v,t,f)^2\big\|_{H_x^{\e}\h^{\ktwo,\ttwo}}^2\Big]
     \\
        &\leq C
        \E\bigg[\!\int_{Q}\int_{\R^{\dirnum}}\!\!\!
        \bigg(\sum_{j=0}^{\e}\sum_{h=0}^{\ktwo}\left|D_y^jD_v^{h+2}\psi\right|^2\bigg) \bigg(\sum_{j=0}^{\e}\sum_{h=0}^{\ktwo}\left|D_y^jD_v^{h}(\sigma(y,v,t,f_{MN})^2-\sigma(y,v,t,f))^2\right|^2\bigg)(1+|v|)^{-2\ttwo}dy\,dv\bigg]
        \\
        &\leq C
        \E\bigg[\!\int_{Q}\int_{\R^{\dirnum}}\!\!\!
        \bigg(\sum_{j=0}^{\e}\sum_{h=0}^{\ktwo}\left|D_y^jD_v^{h+2}\psi\right|^2\bigg)\big(\w_1(Q\times\R^{\dirnum})(f_{MN}(t),f(t))\big)^2 (1+|v|)^{2-2\ttwo}dy\,dv\bigg]
        \\
        &\leq C
        \|\psi\|_{H_x^{\e}\h^{\ktwo+2,\ttwo-1}}\E\Big[\big(\w_2(Q\times\R^{\dirnum})(f_{MN}(t),f(t))\big)^2\Big].
    \end{split}
\end{align}
}
Now Theorem \ref{theorem/ convergence of the single and joint empirical measure} ensures that the last line vanishes uniformly in $t\in[0,T]$  as $M,N\to\infty$.
\end{proof}


\section*{Acknowledgements}
This research has been supported by the EPSRC Centre for Doctoral Training in Mathematics of Random Systems: Analysis, Modelling and Simulation (EP/S023925/1) and by the Advanced Grant Nonlocal-CPD (Nonlocal PDEs for Complex Particle Dynamics: Phase Transitions, Patterns and Synchronization) of the European Research Council Executive Agency (ERC) under the European Union's Horizon 2020 research and innovation programme (grant agreement No. 883363).
The author deeply thanks José A. Carrillo, Pierre Roux and Benjamin Fehrman for comments and helpful discussions.

\bibliographystyle{abbrv}

\end{document}